\documentclass[11pt]{article}
\usepackage[utf8]{inputenc}

\usepackage{comment}
\usepackage[numbers,sort&compress]{natbib}

\usepackage{geometry}
\usepackage{amsthm,amssymb,amsmath,mathtools}

\usepackage[spaces,hyphens]{url}

\usepackage{hyperref}
\usepackage{blindtext}
\usepackage{titlesec}
\usepackage{bbm}
\hypersetup{
  colorlinks   = true, %Colours links instead of ugly boxes
  urlcolor     = blue, %Colour for external hyperlinks
  linkcolor    = black, %Colour of internal links
  citecolor   = red %Colour of citations
}
\usepackage{pdfpages}
\usepackage[all]{hypcap}
\usepackage[font=small]{caption}
\usepackage[rightcaption]{sidecap}
\theoremstyle{plain}
\newtheorem{theorem}{Theorem}[section]
\newtheorem{lemma}[theorem]{Lemma}
\newtheorem{corollary}[theorem]{Corollary}
\newtheorem{proposition}[theorem]{Proposition}
\theoremstyle{definition}
\newtheorem{remark}[theorem]{Remark}

\theoremstyle{definition}

\numberwithin{equation}{section}

\title{Particle systems and the supercooled Stefan problem: Long time behaviour and fluctuations}
\author{Thomas Blore\thanks{Mathematical Institute, University of Oxford. \textit{Email:} thomas.blore@maths.ox.ac.uk },  D.G.M. Flynn\thanks{Mathematical Institute, University of Oxford. \textit{Email:} donald.flynn@maths.ox.ac.uk } and Ben Hambly\thanks{Mathematical Institute, University of Oxford. \textit{Email:} ben.hambly@maths.ox.ac.uk }} 

\begin{document}
\maketitle

\begin{abstract}
    We consider an infinite system of particles on the positive real line, initiated from a Poisson point process, which move according to Brownian motion up until the hitting time of a barrier. The barrier increases when it is hit, allowing for the possibility of sequences of successive jumps to occur instantaneously. The scaling limit gives a representation for the supercooled Stefan problem and its free boundary, from which we can determine the asymptotic limit for the barrier movement. We establish also a functional central limit theorem describing the fluctuations about this limit, from which we can determine the rate of convergence for particle system approximations to the supercooled Stefan problem under ``weak feedback" assumptions introduced in \cite{5}. We finally show that these fluctuations are identical to the fluctuations of a continuous analogue to 1D-MDLA.
\end{abstract}

%\tableofcontents
\section{Introduction}

In this work, we investigate a model for aggregation, given by the following particle system. Let $a>0$ be a positive constant. At time zero, a random configuration of particles is placed via a Poisson point process of unit intensity on the positive half line at points
$$0\leq W_0^{1}\leq \ldots \leq W_0^{n}\leq \ldots .$$
Particles are diffusive in the sense that, independently of the initial particle configuration, the time $t$ position of particle $i$ is given by
$$W^{i}_t=W^{i}_0+B^{i}_t,$$
where each $B^{i}$ is an independent standard Brownian motion. 

\begin{figure}
    \centering
    \includegraphics[width=0.9\linewidth]{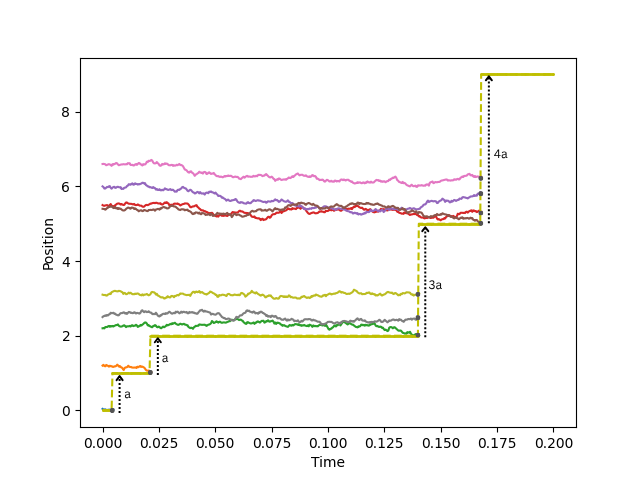}
    \caption{Dynamics for the particle system, and moving barrier, with $a=1$ and the lowest $9$ initial particles.}
    \label{fig1}
\end{figure}

The system of particles is coupled to a moving barrier with position $\xi_t$ at time $t$ that represents the front of a growing aggregate. This is defined recursively in the following manner. At time $0$, we set $\xi_0=0$. $\xi_t$ remains at 0, until time $T_0$, the first hitting time of 0 by any particle. At this time, the barrier height increases to $ak$, where $k$ is the smallest integer such that the interval $[\xi_{T_0^-},\xi_{T_0^-}+ak]$ contains $k$ particles. This captures the idea that when a particle hits the barrier, it jumps by $a$ and if the barrier moves over more particles, they each cause a jump of $a$. Any particles that are now below $\xi_t$ are referred to as ``absorbed" and have no further effect on the movement of $\xi$. The barrier remains at $\xi_{T_0}$ until $T_1$, the first hitting time of $\xi_{T_0}$ by any non-absorbed particle. At time $T_1$, the barrier height again increases by $ak$, where $k$ is the smallest integer such that the interval $[\xi_{T_1^-},\xi_{T_1^-}+ak]$ contains $k$ particles, and these particles are ``absorbed". This mechanism for barrier growth is repeated at each subsequent hitting time $T_n$ of $\xi_{T_{n-1}}$ by any non-absorbed particle, which leads to the following equations for the system dynamics:
\begin{align}
\label{xidefnequation}
    &W_t^i=W_0^i+B_t^i,\nonumber\\
    &\tau_i=\inf\{t\geq 0:W_t^i\leq \xi_t\},\nonumber\\
    &\xi_t=a\sum_{i=1}^\infty \mathbbm{1}_{\{\tau_i\leq t\}},\\
    &\Delta \xi_t=a\inf\{x:\rho_{t^-}([\xi_{t^-},\xi_{t^-}+a x])<x\},\nonumber\\
    &\rho_{t^-}= \left(\sum_{i=1}^\infty\delta_{W_{t}^i}\mathbbm{1}_{\{\tau_i\geq t\}}\right).\nonumber
\end{align}
A simulation of the system is shown in Figure \ref{fig1}.

In \cite{7} this model is studied for $a<1$. The authors of \cite{7} conjecture that under a rescaling, the barrier height
$\xi_t/\sqrt{t}$ converges to some distribution. We establish a strong form of this conjecture through a connection to the supercooled Stefan problem, a simplified model for the freezing of a supercooled liquid, which has recently attracted attention in the probability literature \cite{2,3,minimal,PDEresult,NSsurfacetension,nadtochiy2022tension}. By adapting a monotonicity property from \cite{Dembo_2019}, and extending previous work on probabilistic formulations of the supercooled Stefan problem \cite{delarue2015particle, 5} to the case of non-integrable initial data, we are able to determine the exact form of the limit for $\xi_t/\sqrt{t}$ when $a<1$, and the fluctuations around this limit. 
The probabilistic representations of the supercooled Stefan problem that we use to investigate the behaviour of $\xi_t$ have previously been applied to the setting of integrable initial data \cite{3, mustapha2023wellposednesssupercooledstefanproblem, minimal}.
\subsection{One dimensional MDLA}
\label{1dmdlaintro}
A closely related system to \eqref{xidefnequation} is given by one-dimensional multi-particle diffusion limited aggregation (1D-MDLA). This is another model for the growth of a sticky aggregate, which has previously been studied in discrete space \cite{elboim2020criticalonedimensionalmultiparticledla,sid,sly2016onedimensionalmultiparticlediffusionlimited}. In continuous space (off-lattice), 1D-MDLA is given by a system of particles $W_t^i$, with identical dynamics to \eqref{xidefnequation}, and which are coupled to a moving barrier with position $\bar \xi_t$ at time $t$. The barrier $\bar \xi_t$ is defined recursively as follows. We set $\bar\xi_0=0$. $\bar\xi_t$ remains at 0, until time $T_0$, the first hitting time of 0 by any particle. At this time, the barrier height increases to $a$. Any particles that are now below $\bar \xi_t$ are again referred to as ``absorbed" and have no further effect on the movement of $\bar \xi$. The barrier remains at $\bar \xi_{T_0}$ until $T_1$, the first hitting time of $\bar\xi_{T_0}$ by any non-absorbed particle. At time $T_1$, the barrier height again increases by $a$. This mechanism for barrier growth is repeated at each subsequent hitting time $T_n$ of $\bar \xi_{T_{n-1}}$ by any non-absorbed particle, which leads to the following equations for the system dynamics: 
\begin{align}
\label{xidladefnequation}
    &W_t^i=W_0^i+B_t^i,\nonumber\\
    &\bar \tau_i=\inf\{t\geq 0:W_t^i\leq \bar \xi_t\},\\
    &\bar \xi_t=a\sum_{i=1}^\infty \mathbbm{1}_{\{\bar \tau_i\leq t,\space W^i_{\bar \tau_i}=\bar \xi_{\bar \tau_i-}\}}.\nonumber
\end{align}

This system is almost identical to \eqref{xidefnequation}, except that when $\bar \xi$ increases, multiple particles could be ``absorbed", without $\bar \xi_t$ further increasing. Due to the similarities between \eqref{xidladefnequation} and \eqref{xidefnequation}, we extend our analysis of the behaviour of $\xi_t/{\sqrt{t}}$, to the behaviour of $\bar \xi_t/{\sqrt{t}}$, and determine that $\bar \xi_t/{\sqrt{t}}$ has an identical limit and identical fluctuations when $a<1$.

\subsection{The supercooled Stefan problem}

The one-dimensional one-phase supercooled Stefan problem can be formulated as the following free boundary problem for a function $V:\mathbb{R}_+\times\mathbb{R}_+\to \mathbb{R}$ and a free boundary $\Lambda:\mathbb{R}_+\to\mathbb{R}_+$:
\begin{align}
    \partial_t V(t,x) &=\frac{1}{2}\partial_{xx} V(t,x), \quad x>\Lambda_t,\nonumber\\
     \dot{\Lambda}_t &=\frac{1}{2}\partial_x V(t,\Lambda_t),\label{pde}\\
    V(t,\Lambda_t) &=0, \qquad\qquad t\geq 0 &\nonumber\\
    V(0,x) &=g(x),\qquad\;\; x>0.\nonumber
\end{align}
The function $-V(t,x)$ is a rescaling of temperature at a time $t$ and position $x$, where the freezing temperature is taken as zero. The free boundary $\Lambda_t$ is the position of the interface between the solid and liquid phases. The initial data $g$ is non-negative, hence the temperature $-V(0,x)\leq 0$ is below the freezing point, leading to growth of the solid phase.  

The limiting behaviour of the barrier $\xi_t$ in our particle system will be determined by the behaviour of the free boundary in the Stefan problem \eqref{pde} for non-integrable initial data $g$. Explicit similarity solutions are known with initial data $g\equiv a$ for $a<1$. The free boundary is given by $\Lambda_t=K_{a}\sqrt{t}$ where the constant $K_a$ is the unique solution to
\begin{align}
\label{kalpha}
a=K_a\sqrt{2\pi} (1-\Phi(K_a))e^{K_a^2/2}.
\end{align} 
We use $\Phi$ to denote the standard Gaussian cumulative distribution function.
\subsection{Probabilistic reformulation of the Stefan problem}

For a generic function $g$, classical solutions to \eqref{pde} need not exist. In particular, the PDE may exhibit blow-up in finite time, where the derivative of $\Lambda$ ceases to exist, and jump discontinuities in the movement of the free boundary may occur \cite{2}. This causes weak-formulations of the system to be ill-posed. This issue was resolved for a large class of initial densities by Delarue et al. \cite{3} via a reformulation of the system as the solution to the following McKean-Vlasov type equation:
\begin{align}
    X_t &=X_{0^-}+B_t-\Lambda_t,\nonumber\\
    \tau &=\inf\{t\geq0:X_t\leq 0\},\label{mv}\\
    \Lambda_t &=\alpha\mathbb{P}(\tau\leq t),\nonumber\\
    \Delta\Lambda_t&=\inf\{x>0:\alpha\mathbb{P}(X_t\in [0,x],\tau\geq t)<x\}.\nonumber
\end{align}
In this equation $\alpha=||g||_{L^1(\mathbb{R^+})}<\infty$, the initial condition $X_{0^-}$ has the probability density $g/\alpha$, and $X_{0^-}$ is independent of $B$, a standard Brownian motion on $\mathbb{R}$. The last equation in \eqref{mv} is called the physical jump condition and captures the jump discontinuities in the free boundary. This provides a weak formulation of the supercooled Stefan problem \eqref{pde}, where the free boundary is given by $\Lambda_t$, and $V(t,x)$ corresponds to the density of the process $X_t+\Lambda_t$ at time $t$, see \cite[Theorem 1.1]{DIRT1}.
%$$V(t,x)\textrm{d}x=\alpha\mathbb{P}(X_t+\Lambda_t\in \textrm{d}x,\tau>t).$$

The probabilistic formulation \eqref{mv} requires $g\in L^1(\mathbb{R}^+)$ and therefore cannot be directly compared to the system \eqref{xidefnequation} which has an infinite mass of particles. We therefore consider a generalisation of the McKean-Vlasov system, given in \cite[Equation 1.2]{laplace}:
\begin{align}
    X^x_t &=x+B_t-\Lambda_t,\nonumber\\
    \tau^x &=\inf\{t\geq 0:X^x_t\leq 0\},\label{alternativetomvr}\\
    \Lambda_t &=\int_0^\infty \mathbb{P}(\tau^x\leq t)g(x)\textrm{d}x.\nonumber
    %\\
    %\tilde V(t,x)\textrm{d}x&=\int_0^\infty\mathbb{P}(y+B_t-\Lambda_t\in \textrm{d}x,\tau^y>t)g(y)dy.\nonumber
\end{align}
As in \eqref{mv}, $B$ is a standard Brownian motion, and $g$ is a positive, measurable function. In general, as in \cite{3} it is also necessary to enforce that the free boundary has ``physical" jumps:
\begin{align}
\label{alternativephysical}
 \Delta \Lambda_t=\inf\left\{x>0:\int_0^\infty\mathbb{P}(X_{t^-}^y\in[0,x],\tau^y\geq t)g(y)\textrm{d}y<x\right\}.
\end{align}

\subsection{Particle system reformulation}
\label{reformulationsec}
In order to relate the behaviour of the McKean-Vlasov system \eqref{alternativetomvr}, to the particle system studied in \eqref{xidefnequation}, we introduce the following family of particle systems which are parametrised by integers $N$:
\begin{align}
 &W_t^{i,N} = W_{0}^{i,N}+B_t^{i,N},\nonumber \\
&{} \tau_{i}^N=\inf \left\{t \geq 0: W_t^{i,N} \leq \Lambda^N_t\right\},
\label{aligned:Particlealphanocentre} \nonumber\\
&\Lambda^N_t  =\frac{1 }{N}\sum_{i=1}^\infty \mathbbm{1}_{\left\{ \tau_{i}^N\leq t\right\}}{}_{,}\\
&\Delta \Lambda^N_t=\frac{1}{N}\inf\left\{k\in\mathbb{N}:\rho^N_{t^-}\left(\left[\Lambda^N_{t^-},\Lambda^N_{t^-}+\frac{k}{N}\right]\right)\leq \frac{k}{N}\right\}_{,}\nonumber\\
&\rho^N_{t^-}=\frac{1}{N}\sum_{i=1}^\infty \delta_{W_t^{i,N}}\mathbbm{1}_{\{t\leq \tau_{i}^N\}}\nonumber.
\end{align}
The points $(W_0^{i,N})_{i\in\mathbb{N}}$ are the points of a Poisson process of intensity $Ng$, which are independent of the collection of i.i.d. standard Brownian motions $(B_t^{i,N})_{i\in\mathbb{N}}$. The measure $\rho^N_{t^-}$ is the empirical measure of the non-absorbed particles just before time $t$. 

In the case that $g\equiv a$, the system \eqref{aligned:Particlealphanocentre} corresponds to a rescaling of \eqref{xidefnequation}. We set $W_t^i=Na W_{t/N^2a^2}^{i,N}$ and $\xi_t=Na\Lambda^N_{t/N^2a^2}$ to obtain a system with identical dynamics to \eqref{xidefnequation}. From this scaling relation, $\xi_{N^2a^2t}=Na\Lambda^N_{t}$, we may transfer information on the long time behaviour of system \eqref{xidefnequation}, with intensity $g\equiv 1$, to information on the system \eqref{aligned:Particlealphanocentre}, and vice versa.

This rescaling may also be applied to the system \eqref{xidladefnequation}, by setting 
\begin{align}
 &W_t^{i,N} = W_{0}^{i,N}+B_t^{i,N},\nonumber \\
&{} \bar \tau_{i}^N=\inf \left\{t \geq 0: W_t^{i,N} \leq \bar \Lambda^N_t\right\},
\label{aligned:Particlealphanocentredla} \\
&\bar \Lambda^N_t  =\frac{1 }{N}\sum_{i=1}^\infty \mathbbm{1}_{\left\{ \bar \tau_{i}^N\leq t,\space W_{\bar \tau_i^N}^{i,N}=\bar \Lambda^N_{\bar \tau_i^N}\right\}}{}_{.}\nonumber
\end{align}
When $g\equiv a$, the process $\bar \xi_t=Na\bar \Lambda^N_{t/N^2a^2}$ has identical dynamics to those given in \eqref{xidladefnequation}.

\subsection{The weak feedback regime}

We work in the weak feedback regime equivalent to that given in \cite{5}. The model satisfies the weak feedback condition if the initial condition, a non-negative measurable function $g$, satisfies the following:
\begin{itemize}
    \item  \hypertarget{W}{(W)} $||g||_\infty<1$.
\end{itemize} We impose this condition on $g$ throughout.

The weak feedback condition considerably simplifies the behaviour of solutions to \eqref{alternativetomvr} as it ensures that the feedback effect from absorption is small. Suppose that $g$ satisfies the weak feedback condition. Then, from the results of \cite[Lemma 2.1]{2}, \cite[Proposition 2.2]{3}, \cite[Lemma 2.1, Theorem 2.2]{5}, and \cite[Theorem 4.1]{ledger2024supercooledstefanproblemtransport}, we observe that a unique solution to \eqref{alternativetomvr} exists and is continuous, with $V(t,x)\leq ||g||_\infty<1$ for all times $t$. Further, from a simple adaptation to the arguments within \cite[Theorem 1.1]{PDEresult}, given in the appendix as Proposition \ref{smoothness}, it can be seen that there is a unique classical solution to the supercooled Stefan problem on the time interval $(0,\infty)$. This corresponds exactly to the unique solution to \eqref{alternativetomvr}.

\subsection{Main Results}

For $g$ which satisfies the weak feedback condition, we bound the deviation of the particle system boundaries given in \eqref{aligned:Particlealphanocentre}, from the free boundary in the supercooled Stefan problem \eqref{pde}. In Theorem \ref{bdtheorem}, we give a sub-Gaussian tail bound on the probability that $\sup_{t\leq T}|\Lambda^N_t-\Lambda_t|>x/\sqrt{N}$. Using this bound, and a rewriting of the system applicable due to the fact that the initial particle configuration is given by a Poisson process, in Theorem \ref{clt} we establish a functional central limit theorem result, demonstrating that the scaled deviation converges to a Gaussian process, solving a Volterra integral equation. 

To obtain the limiting behaviour of the system \eqref{xidefnequation} we choose $g\equiv a$ with $a<1$. Thus \eqref{xidefnequation} satisfies the weak feedback condition, and using the rescaling from \eqref{aligned:Particlealphanocentre}, we have the following version of Theorem~\ref{clt} for $\xi_t/\sqrt{t}$ and 
%the fluctuations around this limit.
its fluctuations.

\begin{theorem}
\label{conjtheorem}
For $a<1$, the barrier $\xi_t$ satisfies
 \[ \xi_{st}/\sqrt{t} \to K_a\sqrt{s}, \mbox{uniformly on compact sets in probability as $t\to\infty$}.  \]
 The constant $K_a$ is given by \eqref{kalpha}. For any $T>0$, there exists a constant $C(a,T)>0$ such that for $x<t^{1/4}$, $$\mathbb{P}(\sup_{s\leq T}|\xi_{st}/\sqrt{t}-K_a\sqrt{s}|>xt^{-1/4})\leq 2e^{-C(a,T) x^2}.$$

Further, for any $T\geq 0$, the process $t^{1/4}(\xi_{st}/\sqrt{t}-K_a\sqrt{s})_{s\leq T}$ converges in distribution, under the uniform topology, to a mean zero Gaussian process in $s$, as $t\rightarrow\infty$. The covariance of this process is given by $\mathbb{K}(s,l)$ where $\mathbb{K}(s,l)$ is the unique solution to 
%$$\mathbb{K}(s,l)=a^2\int_0^s\int_0^l\mathbb{K}(z,y)f_{\tau_s}(z)f_{\tau_l}(y)dz\textrm{d}y+K_a\sqrt{s\wedge l},$$
\begin{align*}\mathbb{K}(s,l)&=a\int_0^s\mathbb{K}(z,l)f_{\tau_s}(z)\textrm{d}z+\mathbb{L}(l,s),\\
\mathbb{L}(s,l)&=a\int_0^s\mathbb{L}(z,l)f_{\tau_s}(z)\textrm{d}z+K_a\sqrt{s\wedge l},\end{align*}
and $f_{\tau_s}$ is the probability density function of the time at which $(W_z+K_a\sqrt{z})_{z\leq s}$ is maximised for a Brownian motion $W$.
\end{theorem}
We also show that $\bar\xi_t/\sqrt{t}$ has the same long time behaviour and fluctuations as $\xi_t/\sqrt{t}$.
\begin{theorem}
\label{dlafluc1}
Suppose that $a<1$. Then
    $$\frac{\sup_{s\leq t}|\bar \xi_{s}-\xi_{s}|}{t^{1/4}}\rightarrow 0\textrm{  in probability, as }t\rightarrow\infty.$$
\end{theorem}
Finally, by a straightforward extension to arguments of Theorem \ref{clt}, we establish a corresponding result for the fluctuations of the rescaled empirical measure. In the limit, the fluctuations solve an SPDE with Gaussian noise, specified in Theorem \ref{fclthm}.

\subsection{Related literature and outline}

The particle system \eqref{xidefnequation} may be considered a continuous space (off-lattice) analogue to the ``frictionless" MDLA studied in \cite{Dembo_2019}. In that model, $N_i$  particles are initially placed at each integer $i$, where the $N_i$ are independent Poisson random variables of mean 1. Particles move as independent, continuous time, simple symmetric random walks on $\mathbb{Z}$; and cause a barrier of height $\xi_t$ to grow when moving into contact with this barrier. The authors consider the critical case where $a=1$, and obtain a scaling limit for $\xi_t/t^{2/3}$. An identical scaling limit for the model \eqref{xidefnequation} with $a=1$ is established in \cite{blore2025particlesystemssupercooledstefan}. 

Such MDLA processes have previously been connected to the supercooled Stefan problem. In the work of \cite{dlachayesswindle}, a version of MDLA in one dimension is considered, where the particles move as a continuous time exclusion process on the lattice. It is then shown that a Stefan problem is obtained in the hydrodynamic limit. Such connections between MDLA and Stefan problems are also shown to hold in dimension two in \cite{dlaNS}, in which the authors prove that MDLA converges to a probabilistic analogue to the Stefan problem. 

The McKean-Vlasov formulation of the supercooled Stefan problem \eqref{mv} arose from problems in neuroscience \cite{DIRT1, delarue2015particle} and in systemic risk \cite{2}. The initial work on the neuroscience variant of this system, in \cite{DIRT1, delarue2015particle}, showed that there could be discontinuities and established existence of solutions under the physical jump condition. Subsequently, further properties of this system have been investigated, such as regularity \cite{2,3}, uniqueness of solutions \cite{5, mustapha2023wellposednesssupercooledstefanproblem, 3}, and different solution concepts \cite{minimal,baker}. 

The supercooled Stefan problem with a cumulative distribution function for its initial condition also arises in a model for cell polarization \cite{Calvez2012, BHJ2025}.

Finally we outline our approach. In Section \ref{deterministic barriers}, we determine how the free boundary in \eqref{aligned:Particlealphanocentre} may be estimated by considering the system with a deterministic boundary. In Section \ref{uniquenesssection}, we prove quantitative estimates on the difference between the free boundary in the particle systems \eqref{aligned:Particlealphanocentre} and its limit. We then prove our main fluctuation result in Section \ref{fclt}, and extend this to 1D-MDLA in Section \ref{mdlasec}. Finally, in Section \ref{fffclt} we establish a functional central limit theorem for the empirical measure. The appendix contains a number of proofs deferred from the main text.

\section{Deterministic barriers}
\label{deterministic barriers}
To determine bounds upon the movement of the barrier, we consider the behaviour of the system with a deterministic barrier $f_t$. If the number of particles absorbed by a barrier of height $f_t$ would be too few to increase the barrier height to the value $f_t$, then $\Lambda^N_t$ must be bounded above by $f_t$.

We define a random operator $\Gamma_N:D([0,\infty),\mathbb{R})\rightarrow D([0,\infty),\mathbb{R})$. For a c\'adl\'ag function $f$ we define $\Gamma_N(f)$ via,
\begin{align*}
    &X_t^{i,N}(f)=W_0^{i,N}+B_t^{i,N}-f_t,\\
    &\tau_{i}^N(f)=\inf\{t\geq 0: X_t^{i,N}(f)\leq 0\},\\
    &\Gamma_N(f)_t:=\frac{1}{N}\sum_{i=1}^\infty\mathbbm{1}_{\{\tau_{i}^N(f)\leq t\}}.
\end{align*}
%In subsequent sections, we apply this operator to linear functions which we denote by $l_t^{c,d}=ct+d$ and $l_t^c=l_t^{c,0}.$

\begin{lemma}
\label{Alemmatoboundathing}
    Consider a fixed $T>0,$ and $f$ a c\'adl\'ag function.  Suppose that $\Gamma_N(f)_t\leq f_t$ for all $t\in[0,T]$, and $f(0)\geq 0$. Then $\Lambda^{N}_t\leq f_t$ for all $t\in[0,T]$. Suppose instead that $\Gamma_N(f)_{t}\geq f_t$ for all $t\in[0,T]$, $\Gamma_N(f_{.-}+x1_{\{.\geq t\}})_t\geq f_t+x$ for all $0\leq x\leq \Delta f_t$ at any time $t$ at which $f$ jumps upwards, and $f(0)\leq 0.$ Then $\Lambda^N_t\geq f_t$ for all $t\in[0,T]$.
\end{lemma}
This result follows through minor modifications to the arguments of \cite[Proposition 2.1]{Dembo_2019}. From this bound it is simple to check that a global in time solution to the particle system \eqref{aligned:Particlealphanocentre} exists, for example by using the arguments of Section \ref{weakproofsection}.
%\begin{remark}
%The upper bound in the above lemma also holds for c\'adl\'ag functions, as well as continuous functions. The lower bound holds for c\'adl\'ag functions, as long as the jumps are only downwards or if at any upwards jump of $f$, $\Gamma_N(f_{.-}+x1_{\{.\geq t\}})_t\geq f_t+x$ for all $x\leq \Delta f_t$.
%\end{remark}

 We define also an averaged operator $\Gamma$ on c\'adl\'ag functions $f$ via
\begin{align}
    &X_t^{i}(f)=W_0^i+B_t^i-f_t,\nonumber\\
    &\tau_i(f)=\inf\{t\geq 0:X_t^{i}(f)\leq 0\},\label{gammadefn}\\
    &\Gamma(f)_t:=\sum_{i=1}^\infty \mathbb{P}(\tau_i(f)\leq t).\nonumber
\end{align}
This is identical to setting $\Gamma(f)_t=\mathbb{E}(\Gamma_1(f)_t)$. An analogous operator was first defined in \cite{DIRT1} and used in \cite{2, NS19, minimal}, and can be used to construct solutions to \eqref{alternativetomvr} by iteratively applying $\Gamma$ to the constant function 0, noting that fixed points of $\Gamma$ correspond to solutions of \eqref{alternativetomvr}. Since this operator is an averaged version of $\Gamma_N$ we can easily bound $\Gamma_N(f)-\Gamma(f)$ and then transfer bounds on $\Gamma(f)$ to bounds on $\Gamma_N(f).$

\section{Quantitative bounds}
\label{uniquenesssection}
The main result of this section is the following quantitative bound:
\begin{theorem}
\label{bdtheorem}
   Let $\Lambda$ be the unique solution to \eqref{alternativetomvr}. Then, for each $T\geq 0$ there exists a constant $K$, such that
    $$\mathbb{P}(\sup_{t\leq T}\sqrt{N}|\Lambda^N_t-\Lambda_t|\geq y)\leq 2e^{-Ky^2/(1+y/\sqrt{N})}.$$
\end{theorem}

As a corollary, we obtain the first part of Theorem \ref{conjtheorem}, which we restate here for convenience.
 
\begin{corollary}
    For $y\leq t^{1/4}$:
    $$\mathbb{P}\left(\left|\xi_t/\sqrt{t}-K_a\right|\geq yt^{-1/4}\right)\leq 2e^{-\frac{(1-a)^2y^2}{(3K_{a}+3)}}.$$
\end{corollary}

\begin{proof}
By setting $g\equiv a$, and using the rescaling $\xi_{t}=Na\Lambda^N_{t/N^2a^2}$, it is sufficient to determine the convergence of $\Lambda^N_t$ for the particle system \eqref{aligned:Particlealphanocentre}. Since $a<1$, $g$ satisfies the weak feedback condition, so that the proof of Theorem~\ref{bdtheorem} may be applied. We obtain the lower bound on the constant $K$ in the result of Theorem~\ref{bdtheorem} by inspecting the bounds given in the proof of Theorem~\ref{bdtheorem}.
\end{proof} 
\subsection{Elementary proof of quantitative bounds}
\label{weakproofsection}

We may prove Theorem \ref{bdtheorem} via an elementary argument.
We write $\kappa=1-||g||_\infty.$ For an arbitrary time $t\leq T$, and arbitrary $\epsilon>0$ we bound \begin{align*}
&\Gamma(\Lambda+\epsilon)_t-\Gamma(\Lambda)_t\\&=\int_0^\infty\left( \mathbb{P}\left(\inf_{s\leq t}(x+B_s-\Lambda_s)\leq \epsilon\right)-\mathbb{P}\left(\inf_{s\leq t}(x+B_s-\Lambda_s)\leq 0\right)\right)g(x)\textrm{d}x,\\
&=\mathbb{E}\left(\int_{-\inf_{s\leq t}(B_s-\Lambda_s)}^{-\inf_{s\leq t}(B_s-\Lambda_s)+\epsilon}g(x)\textrm{d}x\right),\\
&\leq \epsilon ||g||_\infty,\\
&=(1-\kappa)\epsilon.
\end{align*}
Since $t\leq T$ was arbitrary, this bound holds for the supremum over $t\leq T$. Further, since $\Gamma(\Lambda)=\Lambda$, we determine that
\begin{align}
\label{boundviakappa}
    \sup_{t\leq T}(\Gamma(\Lambda+\epsilon)_t-\Lambda_t-\epsilon)\leq -\kappa \epsilon.
\end{align}
An identical argument establishes also that
\begin{align}\label{boundviakappalower}\inf_{t\leq T}(\Gamma(\Lambda-\epsilon)_t-\Lambda_t+\epsilon)\geq \kappa\epsilon.\end{align}

Recall the notation $\Gamma_N$ from above Lemma \ref{Alemmatoboundathing}. Consider the event $A^{N,+}_{\epsilon}$ given by
$$A^{N,+}_\epsilon=\left\{\sup_{t\leq T}\left(\Gamma_N(\Lambda+\epsilon)_t-\Gamma(\Lambda+\epsilon)_t\right)\leq \kappa \epsilon\right\}.$$
Applying the bound of \eqref{boundviakappa}, on this event it follows that $$\Gamma_N(\Lambda+\epsilon)_t\leq \Gamma(\Lambda+\epsilon)_t+\kappa\epsilon\leq \Lambda_t+\epsilon \textrm{ for all }t\leq T.$$ By Lemma 
\ref{Alemmatoboundathing} it then follows that $\Lambda^N_t\leq \Lambda_t+\epsilon$ for $t\leq T$. 

This argument can be repeated for $\Lambda-\epsilon$ by considering the event
$$A^{N,-}_\epsilon=\left\{\sup_{t\leq T}\left(\Gamma(\Lambda-\epsilon)_t-\Gamma_N(\Lambda-\epsilon)_t \right)\leq \kappa\epsilon\right\}.$$
Applying the bound from \eqref{boundviakappalower}, on this event it holds that $\Lambda^N_t\geq \Lambda_t-\epsilon$ for $t\leq T$.
We can therefore obtain bounds on the probability that $|\Lambda_t^N-\Lambda_t|>\epsilon$, by bounding the probability of the sets $A_{\epsilon}^{N,+},$ $A_{\epsilon}^{N,-}$.

\begin{lemma} \label{lem:aplus}
    $$\mathbb{P}(A^{N,+}_{y/\sqrt{N}})\geq 1-e^{\frac{-\kappa ^2y^2}{2\Gamma(\Lambda+1)_T+\kappa y/\sqrt{N}}}. $$
\end{lemma}

\begin{proof}
Consider an arbitrary deterministic continuous function $\tilde \Lambda$. 
We define a filtration $(\mathcal{F}_{t})_{t\leq T}$ by $\mathcal{F}_t=\sigma(\Gamma_N(\tilde \Lambda)_s,\space s\leq t)$. We write $\tau_i^N(\tilde \Lambda)=\inf\{t\geq 0:W_0^{i,N}+B_t^{i,N}-\tilde \Lambda _t\leq 0\}$.
We may then check that $\Gamma_N(\tilde \Lambda)_t-\Gamma(\tilde \Lambda)_t$ is a martingale; 
\begin{align*}
    &\mathbb{E}\left(\Gamma_N(\tilde \Lambda)_t-\Gamma(\tilde \Lambda)_t\mid \mathcal{F}_s\right)_,\\
    &\qquad=\mathbb{E}\left(\frac{1}{N}\sum_{i=1}^\infty\mathbbm{1}_{\{s<\tau_i^N(\tilde \Lambda)\leq t\}}\mid \mathcal{F}_s\right)+\Gamma_N(\tilde \Lambda)_s-{\Gamma(\tilde\Lambda)_t}_,\\
    &\qquad=\Gamma(\tilde\Lambda)_t-\Gamma(\tilde\Lambda)_s+\Gamma_N(\tilde\Lambda)_s-{\Gamma(\tilde\Lambda)_t}_,\\
    &\qquad=\Gamma_N(\tilde\Lambda)_s-\Gamma(\tilde\Lambda)_s.
\end{align*}
The second equality above follows since by standard results on the thinning of Poisson processes, the sum $\sum_{i=1}^\infty\mathbbm{1}_{\{s<\tau_i^N(\tilde \Lambda)\leq t\}}$ is a Poisson random variable of mean $N(\Gamma(\tilde\Lambda)_t-\Gamma(\tilde\Lambda)_s)$, which is independent of $\mathcal{F}_s$. 

For arbitrary $\lambda>0$, we apply Doob's inequality to the submartingale \\$\exp\left(\lambda\left(\Gamma_N(\tilde \Lambda)_t-\Gamma(\tilde \Lambda)_t\right)\right)$. By this application, for any $\lambda>0$,
$$\mathbb{P}(\sup_{t\leq T}(\Gamma_N(\tilde \Lambda)_t-\Gamma(\tilde \Lambda)_t)\geq y)\leq \mathbb{E}(e^{\lambda(\Gamma_N(\tilde \Lambda)_T-\Gamma(\tilde \Lambda)_T-y)}).$$
By definition $N\Gamma_N(\tilde\Lambda)_T$ is a Poisson random variable of mean $N\Gamma(\tilde \Lambda)_T$. Using the same arguments as for deriving the Chernoff bound for Poisson random variables, we determine that
$$\mathbb{P}(\sup_{t\leq T}\sqrt{N}(\Gamma_N(\tilde \Lambda)_t-\Gamma(\tilde \Lambda)_t)\geq x)\leq e^{-\frac{x^2}{(2\Gamma(\tilde \Lambda)_T+x/\sqrt{N})}}.$$
Taking $\tilde \Lambda=\Lambda+y/\sqrt{N}$, we have
$$\mathbb{P}(A^{N,+}_{y/\sqrt{N}})\geq 1-e^{-\frac{\kappa^2y^2}{(2\Gamma(\Lambda+1)_T+ \kappa y/\sqrt{N})}}.$$
\end{proof}
Using identical arguments, we can also control the probability of the corresponding lower bound event.
\begin{lemma}\label{lem:aminus}
     $$\mathbb{P}(A^{N,-}_{y/\sqrt{N}})\geq 1-e^{\frac{-\kappa^2 y^2}{2\Gamma(\Lambda+1)_T+\kappa y/\sqrt{N}}}. $$
\end{lemma}
Combining the results of Lemmas~\ref{lem:aplus} and~\ref{lem:aminus}, we may then prove Theorem~\ref{bdtheorem}.
\begin{proof}[Proof of Theorem~\ref{bdtheorem}.]
    From the remarks above Lemma~\ref{lem:aplus}, for arbitrary $y>0$, we have
    $$A^{N,+}_{\frac{y}{\sqrt{N}}}\subset \left\{\sup_{t\leq T}(\Lambda^N_t-\Lambda_t)\leq\frac{y}{\sqrt{N}}\right\} ,$$
    $$A^{N,-}_{\frac{y}{\sqrt{N}}}\subset \left\{-\frac{y}{\sqrt{N}}\leq \inf_{t\leq T}(\Lambda^N_t-\Lambda_t)\right\}.$$
    Consequently:
    $$\mathbb{P}\left(\sup_{t\leq T}\sqrt{N}|\Lambda^N_t-\Lambda_t|>y\right)\leq \mathbb{P}\left(\left(A^{N,+}_{\frac{y}{\sqrt{N}}}\right)^c\right)+\mathbb{P}\left(\left(A^{N,-}_{\frac{y}{\sqrt{N}}}\right)^c\right).$$

    Applying the results of Lemmas~\ref{lem:aplus} and~\ref{lem:aminus}, we then determine that
    $$\mathbb{P}\left(\sup_{t\leq T}\sqrt{N}|\Lambda^N_t-\Lambda_t|>y\right)\leq 2e^{\frac{-\kappa^2 y^2}{\Gamma(\Lambda+1)_T+\kappa y/\sqrt{N}}}.$$
\end{proof}

\begin{remark}
The martingale argument used above requires particles to be initially placed via a Poisson point process. However, for $g\in L^1(\mathbb{R}^+)$, in previous works (e.g.  \cite{delarue2015particle}), exactly $N$ particles are initially placed. The initial positions of these particles are independent and identically distributed random variables which have a distribution corresponding to the density $g$. In this case, the results of this section follow (with possibly different constants) from an application of the Dvoretzky–Kiefer–Wolfowitz inequality.
\end{remark}
\section{Functional central limit theorem}
\label{fclt}

In this section we identify the limiting behaviour of the process $\sqrt{N}(\Lambda^N_t-\Lambda_t)$. 
%Recall that from simple modification of \cite[Theorem 1.1]{PDEresult}, there is a unique classical solution to the supercooled Stefan problem on $(0,T)$ for any $T$. Therefore the result of Remark \ref{remark on smooth case} may be applied to determine that $\sqrt{N}\sup_{t\leq T}|\Lambda^N_t-\Lambda_t|$ is stochastically bounded.
%\begin{remark}
%     The existence of a classical solution to the supercooled Stefan problem is the key assumption used below, while the bound on $g$ is used to simplify the analysis and is stricter than necessary.
%\end{remark}
We use $\delta\Gamma(f,x)$ to denote the directional derivative of the operator $\Gamma$ at $f$, in the direction of the function $x$:
$$\delta\Gamma(f,x)_t=\lim_{\epsilon\rightarrow 0}\frac{\Gamma(f+\epsilon x)_t-\Gamma(f)_t}{\epsilon} \textrm{ for }t\in(0,T).$$

\begin{theorem}
\label{clt}
    Fix $T\geq 0$. The process $\sqrt{N}(\Lambda^N_t-\Lambda_t)_{t\leq T}$ converges in distribution, under the uniform topology, to the unique solution to the Volterra equation
    $$c_t=\int_0^tc_s\mathbb{E}(g(W_s+\Lambda_s)|\tau_t=s)f_{\tau_t}(s)\textrm{d}s+M_t.$$ $W_s$ is a standard Brownian motion, $M_t$ is equal to $B_{\Lambda_t}$ for a standard Brownian motion $B$, and $f_{\tau_t}$ is the probability density function of the time at which the maximum of $W_s+\Lambda_s$ is attained in $[0,t]$. This can equivalently be written as $$c_t=\delta\Gamma(\Lambda,c)_t+M_t.$$
    The solution to this equation is a continuous Gaussian process, started from zero, and with zero mean. 
    
    %The covariance of $x,M$ can be uniquely determined by solving
    
    %by the unique solutions to the following set of equations:
    %\begin{align*}
    %    &\mathbb{E}(M_tx_s)=\int_0^s\mathbb{E}(M_tx_l)\mathbb{E}(g(W_l+\Lambda_l)|\tau_s=l)f_{\tau_s}(l)dl+\Lambda_s \textrm{ for } t\geq s,\\
    %    & \mathbb{E}(M_tx_s)=\int_0^s\mathbb{E}(M_tx_l)\mathbb{E}(g(W_l+\Lambda_l)|\tau_s=l)f_{\tau_s}(l)dl+\Lambda_t \textrm{ for } t<s,\\
    %    &\mathbb{E}(x_tx_s)=\int_0^t\int_0^s\mathbb{E}(x_lx_y))\mathbb{E}(g(W_l+\Lambda_l)|\tau_t=l)f_{\tau_t}(l)\mathbb{E}(g(W_y+\Lambda_y)|\tau_s=y)f_{\tau_s}(y)dl\textrm{d}y+\Lambda_{t\wedge s}.
    %\end{align*}

    %I have just removed the covariance since there doesn't seem to be any insight to be gained.
\end{theorem}

As a corollary, we obtain the second part of Theorem \ref{conjtheorem}, which we restate here for convenience.

\begin{corollary}
\label{4.2}
   Suppose that $a<1$ for the particle system \eqref{xidefnequation}, and fix $T\geq 0$. The process $t^{1/4}(\xi_{st}/\sqrt{t}-K_a\sqrt{s})_{s\leq T}$ converges in distribution, under the uniform topology, to a mean zero Gaussian process in $s$, as $t\rightarrow\infty$. The covariance of this process is given by $\mathbb{K}(s,l)$ where $\mathbb{K}(s,l)$ is the unique solution to 
%$$\mathbb{K}(s,l)=a^2\int_0^s\int_0^l\mathbb{K}(z,y)f_{\tau_s}(z)f_{\tau_l}(y)dz\textrm{d}y+K_a\sqrt{s\wedge l},$$
\begin{align*}\mathbb{K}(s,l)&=a\int_0^s\mathbb{K}(z,l)f_{\tau_s}(z)\textrm{d}z+\mathbb{L}(l,s),\\
\mathbb{L}(s,l)&=a\int_0^s\mathbb{L}(z,l)f_{\tau_s}(z)\textrm{d}z+K_a\sqrt{s\wedge l},\end{align*}
and $f_{\tau_s}$ is the probability density function of the time at which $(W_z+K_a\sqrt{z})_{z\leq s}$ is maximised for a Brownian motion $W$.
\end{corollary}

\begin{proof}
Recall that $\xi_t=Na\Lambda^N_{t/N^2a^2}$. It is sufficient to determine the convergence of $\sqrt{N}(\Lambda^N_{s}-K_a \sqrt{s})$ for the particle system \eqref{aligned:Particlealphanocentre}, with initial density $g(x)=a$ for all $x>0.$ From the result of Theorem \ref{clt}, the convergence of $\sqrt{N}(\Lambda^N_{s}-K_a \sqrt{s})$ to a mean zero Gaussian process holds. The covariance can be explicitly computed using the definition of the limit process given in Theorem \ref{clt}.
\end{proof} 

While there is a large literature on functional limits for McKean-Vlasov particle systems, general results for path dependent interactions require greater regularity than the hitting time map used in defining the process \eqref{aligned:Particlealphanocentre}, hence we give a full proof of the result. %For a short review of such convergence results, see \cite[Section 1]{chaintron2024quasicontinuitymethodmeanfielddiffusions}. 
The conclusions of this section are reliant upon \eqref{alternativetomvr} forming a classical solution to the supercooled Stefan problem, so that the density of particles near zero can be controlled. This holds under the weak feedback condition, as stated in the introduction.

\subsection{A rewriting of the moving barrier system}
\label{rewrite}
To prove Theorem \ref{clt}, we begin with an alternative representation of the moving barrier in the particle system, applicable due to the initial configuration of particles being placed via a Poisson point process. A similar reformulation is applied in \cite{sly2016onedimensionalmultiparticlediffusionlimited} to the discrete space one-dimensional MDLA.

Let $S>0$ be arbitrary, and $(f_t)_{t\in [0,S]}$ be an increasing, c\'adl\'ag  jump process on $[0,S]$. For $t\leq S$, we define densities $\tilde V^f(t,x)$ and $\tilde V^f(t^-,x)$ via $$\tilde V^f(t,x)\textrm{d}x:=\int_0^\infty\mathbb{P}(y+B_t-f_t\in \textrm{d}x,\inf_{s\leq t}(y+B_s-f_s)>0)g(y)\textrm{d}y,$$
$$\tilde V^f(t^-,x)\textrm{d}x:=\int_0^\infty\mathbb{P}(y+B_t-f_{t^-}\in \textrm{d}x,\inf_{s< t}(y+B_s-f_s)>0)g(y)\textrm{d}y.$$
We shall write $\tilde V^f(t,[a,b]),$ $\tilde V^f(t^-,[a,b])$ for the corresponding integral of these densities over the interval $[a,b]$.

Conditioning on $(\Lambda^N_s)_{s\in [0,t)}$, the random measure 
$$\sum_{i=1}^\infty\delta_{(W_{t-}^{i,N}-\Lambda^N_{t-})}\mathbbm{1}_{\{\tau_i^N\geq t\}}$$
is a Poisson random measure of intensity $N\tilde V^{\Lambda^N}(t^-,x)$. Conditional on $(\Lambda^N_s)_{s\in[0,t]}$, the random measure $$\sum_{i=1}^\infty\delta_{(W_t^{i,N}-\Lambda^N_t)}\mathbbm{1}_{\{\tau_i^N>t\}}$$ is a Poisson random measure of intensity given by $N\tilde V^{\Lambda^N}(t,x)$. Finally, conditioning on both the occurrence of a jump at time $t$, and on $(\Lambda^N_s)_{s\in[0,t)}$, the measure $$\sum_{i=1}^\infty\delta_{(W_t^{i,N}-\Lambda^N_{t^-})}\mathbbm{1}_{\{\tau_i^N\geq t,W_t^{i,N}>\Lambda^N_{t^-}\}}$$ is a Poisson random measure, of intensity given by $N\tilde V^{\Lambda^N}(t^-,x)$. This property can be seen iteratively by considering each jump of the process.

Suppose that this property holds up to and including the time of the $i$th jump. It is then clear that the property is preserved between the times of the $i$th and $(i+1)$th jump. At this time, which we label $t$, the jump size $\Delta \Lambda^N_t$ is equal to $$\inf\left\{x:\frac{1}{N}P^t\left(x\right)+1/N-x<0\right\}$$ where $P^t$ is the Poisson process $P^t(x)=\sum_{i=1}^\infty\mathbbm{1}_{\{(W_{t-}^{i,N}-\Lambda^N_{t-})\leq x\}}\mathbbm{1}_{\{\tau_i^N\geq t\}}.$  Therefore, the jump size is a stopping time with respect to the natural filtration of the Poisson process $P^t$. Thus, the point process 
$P^t(x+\Delta \Lambda^N_t)-P^t(\Delta \Lambda^N_t)$ remains conditionally Poisson, with intensity $NV^{\Lambda^N}(t^-,x+\Delta \Lambda^N_t)$.

Applying these formulae, and applying the definition of the minimal jump condition, conditional on the occurrence of a jump at time $t$ and on $(\Lambda_s^N)_{s<t}$, the height of this jump is distributed identically to $$\inf\{x:\frac{1}{N}P\left(N \tilde V^{\Lambda^N}(t^-,[0,x])\right)+1/N-x<0\},$$
where $P$ is a unit rate Poisson process, independent of $(\Lambda_s^N)_{s<t}$. Further, it is simple to see that $(\Lambda^N_t,(\tilde V^{\Lambda^N}(t,.))_{x\in \mathbb{R^+}})$ is a piecewise deterministic Markov process under the natural filtration, and which jumps at rate $-\partial_s \int_0^\infty \tilde V^{\Lambda^N}(s,x)\textrm{d}x|_{s=t^-}=\frac{1}{2}\partial_x\tilde V^{\Lambda^N}(t^-,0)$. 

We can therefore rewrite $\Lambda_t^N$ in terms of Poisson random measures, without referring directly to the particle system. Let $Q$ be a unit intensity Poisson random measure on $(\mathbb{R}^+)^2$, and let $(P_i)_{i\in\mathbb{N}}$ be a collection of independent, unit intensity Poisson processes on $\mathbb{R}^+$, independent of $Q$. We define iteratively a c\'adl\'ag process $(\Lambda^N_t)_{t\geq 0}$ via
\begin{align}
    &\Lambda^N_t:=
    \int_0^t\int_0^{\frac{N}{2}\partial_x\tilde V^{\Lambda^N}(s^-,0)}\Delta \Lambda^N_sQ(\textrm{d}w,\textrm{d}s), \label{newrep}\\
    &\Delta \Lambda^N_s:=\inf\left\{x:\frac{1}{N}P_{N\Lambda^N_{s^-}}\left(N \tilde V^{\Lambda^N}(s^-,[0,x])\right)+1/N-x<0\right \} \textrm{ for }s \textrm{ a jump time.}\nonumber
\end{align}

In the first line above, it can be seen that jumps occur at a rate of $\frac{N\partial_x\tilde V^{\Lambda^N}(s^-,0)}{2}$, and with a height of $\Delta \Lambda_s^N$. The jump size $\Delta\Lambda_s^N$ is defined using only the information of $(\Lambda^N_{l})_{l<s}$ and the Poisson process $P_{N\Lambda^N_{s^-}}$. Therefore, the above process $\Lambda^N_t$ is well defined. This process $(\Lambda^N_t)_{t\geq 0}$ is identical in distribution to the process defined in \eqref{aligned:Particlealphanocentre}. 

\subsection{Consequences of the representation}

We first note that 
\begin{align*}
    \Gamma(\Lambda^N)_t&=\int_0^t\int_0^{\frac{N}{2}\partial_x\tilde V^{\Lambda^N}(s^-,0)}\tilde V^{\Lambda^N}(s^-,[0,\Delta\Lambda_s^N])Q(\textrm{d}w,\textrm{d}s)\\
    & \qquad+\int_0^t\int_0^{\frac{N}{2}\partial_x\tilde V^{\Lambda^N}(s^-,0)}\frac{1}{N}\textrm{d}w\textrm{d}s.
\end{align*}
Recall that the function $\Gamma(f)_t$ returns the amount of mass lost when particles moving according to Brownian motion are killed at position $(f_s)_{s\leq t}$. The first term on the right hand side is then the mass lost when the barrier jumps -  exactly $\tilde V^{\Lambda^N}(s^-,[0,\Delta\Lambda_s^N])$ at each jump time - while the second term gives the mass lost when the barrier remains stationary.

Using the above equality, we can rewrite the representation \eqref{newrep} in the following manner:
\begin{align}
    \Lambda^N_t&=\int_0^t\int_0^{\frac{N}{2}\partial_x\tilde V^{\Lambda^N}(s^-,0)}\left(\Delta \Lambda^N_{s}-\frac{1}{N}\right)Q(\textrm{d}w,\textrm{d}s)
    %\nonumber\\
    %&\qquad 
    +\int_0^t\int_0^{\frac{N}{2}\partial_x\tilde V^{\Lambda^N}(s^-,0)}\frac{1}{N}Q(\textrm{d}w,\textrm{d}s),\nonumber\\
    &=\Gamma(\Lambda^N)_t\nonumber\\
    &\qquad+\int_0^t\int_0^{\frac{N}{2}\partial_x\tilde V^{\Lambda^N}(s^-,0)}\left(\Delta \Lambda^N_{s}-\frac{1}{N}-\tilde V^{\Lambda^N}(s^-,[0,\Delta\Lambda^N_s])\right)Q(\textrm{d}w,\textrm{d}s)\label{rewriteeqtn}\\
    &\qquad +\int_0^t\int_0^{{\frac{N}{2}\partial_x\tilde V^{\Lambda^N}(s^-,0)}}\frac{1}{N}(Q(\textrm{d}w,\textrm{d}s)-\textrm{d}w\textrm{d}s).\nonumber
\end{align}

From this representation we can give a heuristic explanation for the form of the limiting process $c_t$ which appears in Theorem \ref{clt}. We define $c^N$ via $$c_t^N:=\sqrt{N}(\Lambda^N_t-\Lambda_t).$$
The above representation of $\Lambda^N$ shows that 
$$c_t^N=\sqrt{N}(\Lambda^N_t-\Lambda_t)=\sqrt{N}(\Gamma(\Lambda^N)_t-\Gamma(\Lambda)_t)+\hat M_t^N$$
for some process $\hat M_t^N$ which is the sum of the second and third terms in \eqref{rewriteeqtn}. Since the process $\hat M_t^N$ is given by an integral against a compensated Poisson random measure plus a term with jumps of mean zero, it can be seen to be a martingale. The quadratic variation can be computed to be close to $\Lambda_t$. Furthermore, due to the assumptions on the initial density, we can expect that the directional derivative of $\Gamma$ at $\Lambda$ exists, and that
\begin{align*}
\sqrt{N}(\Gamma(\Lambda^N)_t-\Gamma(\Lambda)_t) &=\sqrt{N}\left(\Gamma\left(\Lambda+\frac{c^N}{\sqrt{N}}\right)_t-\Gamma(\Lambda)_t\right),\\
&=\delta\Gamma(\Lambda,c)_t+o(1)_{N\rightarrow\infty}.
\end{align*}

The rest of this section is then dedicated to checking each of these heuristics, by controlling the terms in \eqref{rewriteeqtn}. The results may be proved using simple methods hence we shall only state the results, giving the proofs in Appendix \ref{fcltappen}.

\begin{lemma}
\label{mbound}
The process
$$M_t^N=\sqrt{N}\int_0^t\int_0^{{\frac{N}{2}\partial_x\tilde V^{\Lambda^N}(s^-,0)}}\frac{1}{N}(Q(\textrm{d}w,\textrm{d}s)-\textrm{d}w\textrm{d}s)$$ 
is a martingale of quadratic variation 
$$\Gamma(\Lambda^N)_t-\int_0^t\int_0^{\frac{N}{2}\partial_x\tilde V^{\Lambda^N}(s^-,0)}\tilde V^{\Lambda^N}(s^-,[0,\Delta\Lambda_s^N])Q(\textrm{d}w,\textrm{d}s).$$
    Furthermore, the second term in the quadratic variation, 
$$\int_0^t\int_0^{\frac{N}{2}\partial_x\tilde V^{\Lambda^N}(s^-,0)}\tilde V^{\Lambda^N}(s^-,[0,\Delta\Lambda_s^N])Q(\textrm{d}w,\textrm{d}s)$$ 
converges to zero in probability as $N\to\infty$.
\end{lemma}
\begin{corollary}
\label{mt}
      $M_t^N$ converges uniformly in distribution to a process $M_t$. $M_t$ is a martingale of quadratic variation $\Lambda_t$, started from $0$. In particular, it may be written as $B_{\Lambda_t}$ for a standard Brownian motion $B$.
\end{corollary}

\begin{lemma}
\label{isamg}
    The process
    $$\bar M_t^N=\sqrt{N}\int_0^t\int_0^{\frac{N}{2}\partial_x\tilde V^{\Lambda^N}(s^-,0)}\left(\Delta \Lambda^N_{s}-\frac{1}{N}-\tilde V^{\Lambda^N}(s^-,[0,\Delta\Lambda^N_s])\right)Q(\textrm{d}w,\textrm{d}s)$$ is a martingale.
\end{lemma}
\begin{corollary}
\label{mg2conv}
    The martingale $\bar M_t^N$ converges to zero, uniformly on $[0,T]$, for any $T>0$.
\end{corollary}

Finally, we must determine the limiting distribution of $\sqrt{N}(\Gamma(\Lambda^N)_t-\Gamma(\Lambda)_t)$. For this we require that the limit of $\sqrt{N}(\Lambda^N_t-\Lambda_t)$ is continuous.

\begin{lemma} \label{tightlemma2}
  The process $\sqrt{N}(\Lambda^N_t-\Lambda_t)_{t\in [0,T]}$ is uniformly tight. \\The process $\sqrt{N}(\Gamma(\Lambda^N)_t-\Gamma(\Lambda)_t)_{t\in[0,T]}$ is also uniformly tight.
\end{lemma}

\begin{lemma}
\label{5.7}
For any $(c_t,y_t)_{t\leq T}$ which is a subsequential limit of 
\\$(\sqrt{N}(\Lambda^N_t-\Lambda_t),\sqrt{N}(\Gamma(\Lambda^N)_t-\Gamma(\Lambda)_t))_{t\leq T},$
almost surely $y_t=\delta\Gamma(\Lambda,c)_t$ where $\delta \Gamma(\Lambda,c)$ is the directional derivative of $\Gamma$ at $\Lambda$, in the direction $c$. The directional derivative $\delta\Gamma(\Lambda,c)_t$ may also be written as $\mathbb{E}^W(c_{\tau_t}g(Y_{\tau_t}))$ where $Y_s=W_s+\Lambda_s$ for a Brownian motion $W$, independent of $(\Lambda^N)_{N\in\mathbb{N}}$, $\tau_t$ is the time at which $Y_s$ attains a maximum on $[0,t]$, and $\mathbb{E}^W$ is the expectation taken over the random process $W$.
\end{lemma}

We may now prove the main result of this section.
\begin{proof}[Proof of Theorem \ref{clt}]
    By Corollary \ref{mt}, Corollary \ref{mg2conv}, Lemma \ref{tightlemma2}, and Lemma \ref{5.7} we consider a subsequence on which $(M_t^{N},\bar M_t^{N},\sqrt{N}(\Lambda^N_t-\Lambda_t),\sqrt{N}(\Gamma(\Lambda^N)_t-\Gamma(\Lambda)_t))_{t\leq T}$ converges in distribution to some $(M_t,\bar M_t,c_t,y_t)_{t\leq T}.$
    From the aforementioned results, we can identify that almost surely $\bar M_t=0$, $y_t=\delta\Gamma(\Lambda,c)_t$, and $M_t=B_{\Lambda_t}$ for a Brownian motion $B$. From equation \eqref{rewriteeqtn} we can finally identify that $$c_t=y_t+M_t+\bar M_t=\delta\Gamma(\Lambda,c)_t+B_{\Lambda_t}.$$
    
    We investigate further the solution to this equation. First, the solution to this equation is unique due to the weak feedback assumption. In particular, for two solutions $c^1,c^2$, 
    \begin{align*}
          \sup_{s\leq t}|c^1_s-c^2_s|\leq \mathbb{E}(g(\sup_{s\leq t}(W_s+\Lambda_s)))\sup_{s\leq t}|c_s^1-c_s^2|\leq ||g||_\infty\sup_{s\leq t}|c_s^1-c_s^2|
    \end{align*}
Moreover, it is simple to construct a solution to the Volterra equation, as an $L^2$ limit of Gaussian processes. Indeed consider a sequence of processes
\begin{align*}
    u_0(t)&=M_t,\\
    u_{n+1}(t)&=\int_0^t\mathbb{E}(g(W_{\tau_t+\Lambda_{\tau_t}})|\tau_t=s)f_{\tau_t}(s)u_{n}(s)\textrm{d}s,\\
    x^n_t&=\sum_{i=0}^nu_{i}(t).
\end{align*}
Since $g$ is bounded, and $f_{\tau_t}$ is integrable, it is standard to check that each of the $u_{n}(t)$ processes are $L^2$ limits of Riemann sums and are Gaussian. Further, $||u_n||_{T}=\sup_{s\leq T}|u_{n}(s)|$ must satisfy 
\begin{align*}
    ||u_{n+1}||_T\leq ||g||_\infty||u_n||_T\leq || g||_\infty^n||u_0||_T.
\end{align*}
It is then simple to see that $x^n$ are Gaussian and uniformly Cauchy, with limit $x_t$ which must solve the Volterra equation. Further, 
$\mathbb{E}(||x-x^n||_T^2)\leq \frac{||g||_\infty^n}{1-||g||_\infty}\mathbb{E}(||M||_T^2)$ which converges to zero as $n\rightarrow\infty$. We then conclude that the process $x_t$ is Gaussian, hence the unique solution to the Volterra equation is Gaussian.
\end{proof}
\section{Application to 1D-MDLA}
\label{mdlasec}
Applying the bounds established in the proof of Theorem \ref{conjtheorem}, we show that the long time behaviour and fluctuations of the 1D-MDLA process, $\bar \xi_t$, are identical to those of $\xi_t$. We restate Theorem \ref{dlafluc1}, using the rescaling $Na\bar \Lambda^N_t={\bar \xi_{N^2a^2t}}$ and $Na\Lambda^N_t=\xi_{N^2a^2t}$.

\begin{theorem}
\label{dlafluc2}
Suppose that $g\equiv a$, for some $a<1$. Then
     $\sup_{t\leq 1}\sqrt{N}|\Lambda^N_t-\bar \Lambda^N_t|$ converges in probability to zero as $N\rightarrow\infty$.
\end{theorem}

We prove Theorem \ref{dlafluc2}, by using a reformulation of $\bar \Lambda^N_t$ similar to that used for $\Lambda^N_t$ in Section \ref{rewrite}. By analogous arguments to those within Section \ref{rewrite}, $\bar \Lambda^N_t$ has identical distribution to
$$\bar \Lambda^N_t:=
    \int_0^t\int_0^{\frac{N}{2}\partial_x\tilde V^{\bar \Lambda^N}(s^-,0)}\frac{1}{N}Q(\textrm{d}w,\textrm{d}s).$$
Rather than coupling $\Lambda^N$ and $\bar\Lambda^N$ through the collection $(W_t^{i,N})_{i\in\mathbb{N}}$, we instead couple them through the Poisson point process $Q$. This simplifies comparison between $\Lambda^N$ and $\bar\Lambda^N$, as we need only keep track of the behaviour of $\Lambda^N$ and $\bar \Lambda^N$, rather than the behaviour of the entire system of particles $(W_t^{i,N})_{i\in\mathbb{N}}.$ Further, since the initial distribution is given by $g\equiv a$, we can also rewrite $\tilde V^f(t^-,x)$ given in Section \ref{rewrite} as
 \begin{align}\tilde V^f(t^-,x)=a\mathbb{P}\left(\inf_{0<s\leq t}\left(x+B_s+(f_{t-}-f_{t-s})\right)>0\right).\label{altdefnv}\end{align}
This is increasing in $x$, and is also increasing in $(f_{t-}-f_{t-s})_{0<s\leq t}$. Using the coupling through the point process $Q$, we can then control the behaviour of $\bar \Lambda^N$ in terms of $\Lambda^N$.
\begin{lemma}
\label{lemmaforboundingL}
For all $0< s\leq t$,
    $$\Lambda^N_{t-}-\Lambda^N_{t-s}\geq \bar \Lambda^N_{t-}-\bar \Lambda^N_{t-s}.$$
    Subsequently, for all $t\geq 0$, and all $x\geq 0$
    $$\tilde V^{\Lambda^N}(t^-,x)\geq \tilde V^{\bar \Lambda^N}(t^-,x) .$$
    Further,
    $$\partial_x\tilde V^{\Lambda^N}(t^-,0)\geq \partial_x\tilde V^{\bar \Lambda^N}(t^-,0).$$
\end{lemma}
\begin{proof}
From the observation that $V^{f}(t^-,x)$ is increasing in $(f_{t-}-f_{t-s})_{0<s\leq t}$, the second and third given properties hold immediately upon assumption of the first property. This stated property trivially holds at time zero. Assume that the given property holds up to the $i$th jump time of $(\Lambda^N_t,\bar \Lambda^N_t)$, which we label $T_i$. At this time, by assumption $\partial_x\tilde V^{\Lambda^N}(t^-,0)\geq \partial_x\tilde V^{\bar \Lambda^N}(t^-,0)$. Therefore, if $\bar \Lambda^N$ was to jump at this time, so must $\Lambda^N$. Furthermore, the jump size $\Delta \Lambda^N_{T_i}$ is at least $ 1/N $, hence is at least as large as $\Delta \bar \Lambda^N_{T_i}$. We label $T_{i+1}$ the time of the next jump. For $T_i< t\leq  T_{i+1}$, and $0<s\leq t$,
\begin{align*}
        \Lambda^N_{t-}-\Lambda^N_{t-s}&=\mathbbm{1}_{\{s> t-T_i\}}\left(\Delta \Lambda^N_{T_i}+\Lambda^N_{T_i-}-\Lambda^N_{T_i-s}\right),\\
        &\geq \mathbbm{1}_{\{s>t-T_i\}}\left(\Delta \bar \Lambda^N_{T_i}+\Lambda^N_{T_i-}-\Lambda^N_{T_i-s}\right),\\
        &\geq \mathbbm{1}_{\{s>t-T_i\}}\left(\Delta \bar \Lambda^N_{T_i}+\bar \Lambda^N_{T_i-}-\bar \Lambda^N_{T_i-s}\right),\\
        &=\bar \Lambda^N_{t-}-\bar \Lambda^N_{t-s}
    \end{align*}

The second inequality follows by the inductive assumption. This shows that the given property holds up to the $(i+1)$th jump. By induction, we have therefore shown that the property holds over all time.
\end{proof}

We denote by $\tilde L^N_t$ the mass of ``absorbed" particles, which do not cause growth of $\bar \Lambda^N$. This is defined by
$$\tilde L^N_t:=\int_0^t\int_0^{\frac{N}{2}\partial_x\tilde V^{\bar \Lambda^N}(s^-,0)}\tilde V^{\bar \Lambda^N}(s^-,[0,1/N])Q(\textrm{d}w,\textrm{d}s).$$
Applying the bounds on $c^N_t$ established in Section \ref{fclt}, we shall later show that $\sqrt{N}
\tilde L^N_T$ is small.
\begin{lemma}
\label{lto0}
    $\sqrt{N}\tilde L^N_T$ converges to zero in probability as $N\rightarrow\infty$.
\end{lemma}

The existence of particles which are ``absorbed", but which do not cause an increase in the barrier position, is the key difference between the system \eqref{xidladefnequation} and \eqref{xidefnequation}. The result of Theorem \ref{dlafluc2} therefore follows from simple arguments, once Lemma \ref{lto0} is established.
\begin{proof}[Proof of Theorem \ref{dlafluc2}]
Applying the definition of the operator $\Gamma$, we rewrite $\bar \Lambda^N_t$ as 
$$\bar \Lambda^N_t=\frac{1}{N}\int_0^t\int_0^{\frac{N}{2}\partial_{x}\tilde V^{\bar \Lambda^N}(s^-,0)}(Q(\textrm{d}w,\textrm{d}s)-\textrm{d}w\textrm{d}s)+\Gamma(\bar \Lambda^N)_t-\tilde L^N_t.$$
From this, we first observe that 
\begin{align*}\bar \Lambda^N_t-\Lambda_t&=\Gamma(\bar \Lambda^N)_t-\Gamma(\Lambda)_t\\
&\qquad +\frac{1}{N}\int_0^t\int_0^{\frac{N}{2}\partial_{x}\tilde V^{\bar \Lambda^N}(s^-,0)}(Q(\textrm{d}w,\textrm{d}s)-\textrm{d}w\textrm{d}s)-\tilde L^N_t.\end{align*}

As in Lemma \ref{lem:aplus}, we can then rearrange this to determine that
$$\sup_{t\leq 1}|\bar\Lambda^N_t-\Lambda^{}_t|\leq \frac{1}{1-a}\left(\sup_{t\leq 1}\left|\frac{1}{N}\int_0^t\int_0^{\frac{N}{2}\partial_{x}\tilde V^{\bar \Lambda^N}(s^-,0)}(Q(\textrm{d}w,\textrm{d}s)-\textrm{d}w\textrm{d}s)\right|+\tilde L_1^N\right).$$
The first term is a martingale of quadratic variation bounded by ${\Lambda^N_1}/{N}$. This converges in probability to zero, hence the martingale also converges uniformly in probability to zero. Furthermore, $\tilde L_1^N$ converges to zero by Lemma \ref{lto0}. We have thus established $\bar \Lambda^N_t-\Lambda_t$ converges uniformly in probability to zero. 

From the rewriting of $\bar \Lambda^N_t$ we also observe,
\begin{align*}
    \sqrt{N}(\Lambda^N_t-\bar \Lambda^N_t)&=\sqrt{N}(\Gamma(\Lambda^N)_t-\Gamma(\bar \Lambda^N)_t)\\&\qquad+\frac{1}{\sqrt{N}}\int_0^t\int_{\frac{N}{2}\partial_x\tilde V^{\bar \Lambda^N}(s^-,0)}^{\frac{N}{2}\partial_{x}\tilde V^{\Lambda^N}(s^-,0)}(Q(\textrm{d}w,\textrm{d}s)-\textrm{d}w\textrm{d}s)+\sqrt{N}\tilde L_t^N+\bar M^N_t.
\end{align*}
Rearranging this in the same way as $\bar \Lambda^N-\Lambda$, we determine that
\begin{align*}&\sqrt{N}\sup_{t\leq 1}|\Lambda^N_t-\bar \Lambda^N_t|\\
&\quad\leq \frac{1}{1-a}\left(\sup_{t\leq 1}\left|\frac{1}{\sqrt{N}}\int_0^t\int_{\frac{N}{2}\partial_x\tilde V^{\bar \Lambda^N}(s^-,0)}^{\frac{N}{2}\partial_{x}\tilde V^{\Lambda^N}(s^-,0)}(Q(\textrm{d}w,\textrm{d}s)-\textrm{d}w\textrm{d}s)+\bar M^N_t\right|+\sqrt{N}\tilde L_1^N\right).\end{align*}
From Corollary \ref{mg2conv}, $\sup_{s\leq 1}|\bar M^N_s|$ converges in probability to zero, while from Lemma \ref{lto0}, $\sqrt{N}\tilde L^N_1$ converges in probability to zero. Finally, the term 
$$\frac{1}{\sqrt{N}}\int_0^t\int_{\frac{N}{2}\partial_x\tilde V^{\bar \Lambda^N}(s^-,0)}^{\frac{N}{2}\partial_{x}\tilde V^{\Lambda^N}(s^-,0)}(Q(\textrm{d}w,\textrm{d}s)-\textrm{d}w\textrm{d}s)$$
is a martingale, with quadratic variation bounded above by $\Lambda_t^N-\bar \Lambda^N_t$. Since we have established that both $\Lambda^N_t$ and $\bar \Lambda^N_t$ converge uniformly in probability to $\Lambda_t$, the quadratic variation converges to zero, hence the martingale term must converge uniformly in probability to zero.
 \end{proof}

We split the proof of Lemma \ref{lto0} into two separate sub-lemmas. 
\begin{lemma}
    \label{lto02}
    $\sqrt{N}\tilde L_\epsilon^N$ converges in probability to zero, as $N\rightarrow\infty$, then $\epsilon\downarrow 0$.
\end{lemma}
\begin{proof}
From the definition of $\tilde L^N_t$ and the result of Lemma \ref{lemmaforboundingL}, we can upper bound $\tilde L_\epsilon^N$ above by
\begin{align*}
    &\int_0^\epsilon\int_0^{\frac{N}{2}\partial_x\tilde V^{\Lambda^{N}}(s^-,0)}\tilde V^{\Lambda^{N}}(s^-,[0,1/N])Q(\textrm{d}w,\textrm{d}s)=\sum_{0<s\leq \epsilon,\Delta \Lambda^N_s>0}\tilde V^{\Lambda^N}(s^-,[0,1/N]).
\end{align*}
We upper bound $\tilde V^{\Lambda^N}(t^-,x)$ by using \eqref{altdefnv} to write 
\begin{align*}
    \tilde V^{\Lambda^N}(t^-,x)
    &=a\mathbb{P}\left(\inf_{0<s\leq t}\left(x+B_s-(\Lambda^N_{t-}-\Lambda^N_{t-s})\right)>0\right),\\
    &\leq a\mathbb{P}\left(\inf_{0<s\leq t}\left(x+2N^{-1/2}\sup_{l< t}|c_l^N|+B_s-(\Lambda_{t-}-\Lambda_{t-s})\right)>0\right),\\
    %&=a\int_0^\infty \mathbb{P}\left(y+B_t-\Lambda^N_{t-}\in dx,\inf_{s<t}(y+B_s-\Lambda^N_s)>0\right)dy,\\
    %&\leq a\int_0^\infty \mathbb{P}\left(y+B_t-\Lambda^N_{t-}\in dx,\inf_{s<t}\left(y+B_s-K_a\sqrt{s}+N^{-1/2}\sup_{s<t}|c_s^N|\right)>0\right)dy,\\
    %&\leq a\int_0^\infty \mathbb{P}\left(y+B_t-K_a\sqrt{t}-N^{-1/2}\sup_{s<t}|c_s^N|\in dx,\inf_{s<t}\left(y+B_s-K_a\sqrt{s}+N^{-1/2}\sup_{s<t}|c_s^N|\right)>0\right)dy,\\
    %&= a\int_{N^{-1/2}\sup_{s<t}|c_s^N|}^\infty \mathbb{P}\left(y+B_t-K_a\sqrt{t}-\in d\left(x+2N^{-1/2}\sup_{s<t}|c_s^N|\right),\inf_{s<t}\left(y+B_s-K_a\sqrt{s}\right)>0\right)dy,\\
    %&\leq a\int_{0}^\infty \mathbb{P}\left(y+B_t-K_a\sqrt{t}-\in d\left(x+2N^{-1/2}\sup_{s<t}|c_s^N|\right),\inf_{s<t}\left(y+B_s-K_a\sqrt{s}\right)>0\right)dy,\\
    &=\tilde V^{\Lambda}(t,x+2N^{-1/2}\sup_{s<t}|c_s^N|),\\
    &=\frac{a}{(1-\Phi(K_a))}\left(\Phi\left(K_a+\frac{2\sup_{s<t}|c_s^N|}{\sqrt{Nt}}+\frac{x}{\sqrt{t}}\right)-\Phi(K_a)\right).
\end{align*}
The inequality follows, since $\tilde V^{f}(t^-,x)$ is increasing in $x$ and in $(f_{t-}-f_{t-s})_{0<s\leq t}$. The final equality follows from the form of the similarity solution to the supercooled Stefan problem.
%The first inequality follows as we reduce the height of the free boundary at which particles are killed to $\Lambda_t-N^{-1/2}\sup_{s<t}|c_s^N|$. The second inequality holds since the function
%$$\int_0^\infty \mathbb{P}\left(y+B_t\in dx,\inf_{s<t}\left(y+B_s-K_a\sqrt{s}+N^{-1/2}\sup_{s<t}|c_s^N|\right)>0\right)dy$$ is monotone increasing in $x$.

Fix arbitrary $\epsilon>\delta>0.$ We upper bound 
\begin{align*}
    &\sum_{\delta<t\leq \epsilon,\Delta \Lambda^N_t>0} \tilde V^{\Lambda^N}(t^-,[0,1/N])\\
    &\qquad\leq  \frac{a}{(1-\Phi(K_a))}\sum_{\delta<t\leq \epsilon,\Delta \Lambda^N_t>0}\int_0^{\frac{1}N}\left(\Phi\left(K_a+\frac{x}{\sqrt{t}}+\frac{2\sup_{s<t}|c_s^N|}{\sqrt{Nt}}\right)-\Phi(K_a)\right)dx,\\
    &\qquad \leq \frac{C_a}{N}\sum_{\delta<t\leq \epsilon, \Delta \Lambda^N_t>0}\left(\frac{\sup_{s<t}|c_s^N|}{\sqrt{Nt}}\right)\wedge 1.
\end{align*}
The constant $C_a$ depends only on $a$. 

We increase this upper bound to 
\begin{align}&C_a\int_\delta^\epsilon \left(\left(\frac{\sup_{s<t}|c_s^N|}{\sqrt{Nt}}\right)\wedge 1\right)\textrm{d}\Lambda^N_t
   \nonumber\\&=  C_a\int_\delta^\epsilon \left(\left(\frac{\sup_{s<t}|c_s^N|}{\sqrt{Nt}}\right)\wedge 1\right)\textrm{d}\Lambda_t+C_a\int_{\delta}^\epsilon  \left(\left(\frac{\sup_{s<t}|c_s^N|}{\sqrt{Nt}}\right)\wedge 1\right)\textrm{d}(\Lambda^N-\Lambda)_t\label{bounder} \end{align}

By integrating by parts, the second term in \eqref{bounder},
$$C_a\int_\delta^\epsilon\left( \frac{\sup_{s<t}|c_s^N|}{\sqrt{Nt}}\wedge 1\right)\textrm{d}(\Lambda^N-\Lambda)_t,$$
is bounded above by 
$$C_a \left(\frac{\sup_{s\leq \epsilon}|c_s^N|}{\sqrt{N\epsilon}}(\Lambda^N-\Lambda)_\epsilon+\int_0^{\epsilon}|\Lambda^N_{t-}-\Lambda_{t-}|\textrm{d}TV\left(\left(\frac{\sup_{s\leq t}|c_s^N|}{\sqrt{Nt}}\right)\wedge 1\right)\right).$$
 The second term is bounded above by
\begin{align*}
   &C_a\left(\int_{\delta}^\epsilon\frac{\left( \sup_{s<t}|c_s^N|\right)^2}{2Nt^{3/2}}\mathbbm{1}_{\left\{\frac{\sup_{s<t}|c_s^N|}{\sqrt{Nt}}\leq 1\right\}}\textrm{d}t+\sum_{\delta\leq t\leq \epsilon}\frac{\sup_{s<t}|c_s^N|}{N\sqrt{t}}\mathbbm{1}_{\left\{\frac{\sup_{s<t}|c_s^N|}{\sqrt{Nt}}\leq 1\right\}}\Delta(\sup_{s\leq t}|c_s^N|)\right),\\
   &\leq C_a\left(\int_{\delta}^\epsilon\frac{ \sup_{s<t}|c_s^N|}{2\sqrt{N}t}\textrm{d}t+\frac{\sup_{s\leq \epsilon}|c_s^N|}{\sqrt{N}}\right).
\end{align*}
Plugging this back into equation \eqref{bounder} we have therefore shown that
\begin{align*}
     &\sqrt{N}\sum_{\delta<t\leq \epsilon,\Delta \Lambda^N_t>0} \tilde V^{\Lambda^N}(t^-,[0,1/N])\\
     &\leq C_a\left( \int_\delta^\epsilon \frac{\sup_{s<t}|c_s^N|}{\sqrt{t}}\textrm{d}\Lambda_t+\frac{2\sup_{s\leq \epsilon}|c_s^N|}{\sqrt{\epsilon}}\sup_{t\leq\epsilon}|\Lambda^N_t-\Lambda_t|+\int_{\delta}^\epsilon\frac{ \sup_{s<t}|c_s^N|}{2t}\textrm{d}t+{\sup_{s\leq \epsilon}|c_s^N|}\right),\\
     &\leq C_a\left( \int_\delta^\epsilon \frac{\sup_{s<t}|c_s^N|}{t}\textrm{d}t+\frac{2\sup_{s\leq \epsilon}|c_s^N|}{\sqrt{\epsilon}}\sup_{t\leq\epsilon}|\Lambda^N_t-\Lambda_t|+{\sup_{s\leq \epsilon}|c_s^N|}\right).
\end{align*}
We have increased $C_a$ in the final line, and used that $\Lambda_t=K_a\sqrt{t}$. The second term on the right hand side of this inequality converges in probability to zero as $N\rightarrow\infty$, by the result of Theorem \ref{bdtheorem}. The third term is independent of the choice of $\delta>0$, and converges in probability to zero as $N\rightarrow\infty$, then $\epsilon\downarrow 0$, again by Theorem \ref{bdtheorem}.

To bound the first term, we consider the expectation of $\sup_{s\leq t}|c_s^N|$. From the proof of Lemma \ref{lem:aplus} and \ref{lem:aminus}, we can bound  $$\mathbb{P}\left((1-a)\sup_{s\leq t}|c^N_s|\geq x\right)\leq 2e^{-\frac{x^2}{2K_a\sqrt{t}+3x/((1-a)\sqrt{N})}}.$$ It follows that there exists a constant $\tilde C_a$ such that $\mathbb{E}(\sup_{s\leq t}|c^N_s|)\leq \tilde C_at^{1/4}$ for all $t\geq N^{-2}$. Taking $\delta=N^{-2}$ we then observe that
\begin{align*}
     \mathbb{E}\left( \int_{N^{-2}}^\epsilon \frac{\sup_{s<t}|c_s^N|}{\sqrt{t}}\textrm{d}\Lambda_t\right)
     &\leq \int_0^\epsilon\frac{\tilde C_a K_a}{2t^{3/4}}\textrm{d}t\\
     &\leq \tilde C_a\epsilon^{1/4}.
\end{align*}
By Markov's inequality, it follows that as $N\rightarrow\infty$, then $\epsilon\downarrow 0$, $ \int_{N^{-2}}^\epsilon \frac{\sup_{s<t}|c_s^N|}{\sqrt{t}}\textrm{d}\Lambda_t$ converges in probability to zero. We have therefore determined that $$\sqrt{N}\sum_{N^{-2}<t\leq \epsilon,\Delta \Lambda^N_t>0} \tilde V^{\Lambda^N}(t^-,[0,1/N])\overset{p}{\rightarrow}0$$ as $N\rightarrow\infty$, and then $\epsilon\downarrow 0$.

Finally, we bound 
$$\sqrt{N}\sum_{t\leq N^{-2},\Delta \Lambda^N_t>0} \tilde V^{\Lambda^N}(t^-,[0,1/N])\leq \sqrt{N}\Lambda^N_{N^{-2}}.$$
From the rescaling given in Section \ref{reformulationsec}, this has the same distribution as $\frac{\xi_{a^2}}{a\sqrt{N}}$. This clearly converges to zero in probability as $N\rightarrow\infty$. Therefore, 
$$\sqrt{N}\sum_{t\leq \epsilon,\Delta \Lambda^N_t>0} \tilde V^{\Lambda^N}(t^-,[0,1/N])\overset{p}{\rightarrow} 0$$ as $N\rightarrow\infty$, then $\epsilon\downarrow 0$. Since this was an upper bound on $\sqrt{N}\tilde L^N_\epsilon$ we are done.
\end{proof}
\begin{lemma}
\label{lto01}
    For any fixed $T\geq \epsilon>0$, $\sqrt{N}(\tilde L_T^N-\tilde L_\epsilon^N)$ converges in probability to zero as $N\rightarrow\infty$.
\end{lemma}
%This allows for control of $\bar \Lambda^N-\Lambda$ due to a similar property to Lemma \ref{Alemmatoboundathing}.
%\begin{lemma}
%Let $f$ be a c\'adl\'ag function, with jumps downward.  Suppose that for $t\leq T$ $\Gamma_N(f)_t\geq f_t+L^N_T$. Then $\bar \Lambda^N_t\geq f_t$ for all $t\leq T$.
%\end{lemma}
%\begin{proof}
%Suppose that $\Gamma_N(f)_t\geq f_t+L^N_T$ for all $t\leq T$. 
%We shall show that $\bar \Lambda^N_t\geq f_t$ for all $t\leq T$. At time zero, this property holds trivially. Since $f$ does not jump upwards, at the first time $t$ such that $f_t>\Lambda^N_{t}$
%$$f_t>\bar \Lambda^N_t=\Gamma_N(\bar \Lambda^N)_t-L^N_t\geq \Gamma_N(f)_t-L_T^N.$$ This contradicts the assumption that $\Gamma_N(f)_t\geq f_t+L^N_T$.
%\end{proof}
\begin{proof}
Analogously to the proof of Lemma \ref{lto02}, we upper bound $$\tilde L^N_T-\tilde L^N_{\epsilon}\leq \sum_{\epsilon\leq t\leq T,\Delta \Lambda^N_t>0} \tilde V^{\Lambda^N}(t^-,[0,1/N]).$$
We recall also from the proof of Lemma \ref{lto02} that
$$\tilde V^{\Lambda^N}(t,x)\leq \frac{a}{(1-\Phi(K_a))}\left(\Phi\left(K_a+\frac{2\sup_{s\leq t}|c_s^N|}{\sqrt{Nt}}+\frac{x}{\sqrt{t}}\right)-\Phi(K_a)\right).$$
   Inspecting this expression, there exists a constant $C_\epsilon$, independent of $N$, such that for all $t\geq \epsilon$, this may be further bounded above by $\tilde V^{\Lambda}(t,x)+C_\epsilon N^{-1/2}\sup_{s\leq T}|c_s^N|$.
     
     In between jumps of $\Lambda^N$, $\tilde V^{\Lambda^N}$ evolves via the heat equation with killing at zero. We write $T_i$ for the $i$th jump time of $\Lambda^N$. For $t\in [\epsilon,T]$ we apply the heat equation with killing at zero to the upper bound on $\tilde V^{\Lambda^N}(t,x)$. For $T_{i-1}\leq t<T_i$ we therefore bound
     $$\tilde V^{\Lambda^N}(t,x)\leq \tilde V^{\Lambda}(t,x)+C_\epsilon N^{-1/2}\sup_{s\leq T}|c_s^N|\textrm{erf}\left(\frac{x}{\sqrt{2(t-T_{i-1}\vee \epsilon)}}\right).$$
     It follows that
     \begin{align}\tilde V^{\Lambda^N}(T_i-,[0,1/N])\leq \tilde C_\epsilon \frac{1}{N^2}+ \tilde C_\epsilon N^{-3/2}\sup_{s\leq T}|c_s^N|\textrm{erf}\left(\frac{1}{N\sqrt{2(T_{i}-T_{i-1}\vee \epsilon)}}\right),\label{eqbd1}\end{align}
    for some constant $\tilde C_\epsilon$, and that 
\begin{align}&\frac{N}{2}\int_{T_{i-1}\vee \epsilon}^t\partial_x\tilde V^{\Lambda^N}(s^-,0)\textrm{d}s\nonumber\\
&\qquad \leq \tilde C_\epsilon\left(N(t-T_{i-1}\vee \epsilon)+N^{-1/2}\sup_{s\leq T}|c_s^N|\sqrt{t-T_{i-1}\vee \epsilon}\right).\label{eqbd2}\end{align}
Fix arbitrary $K\geq 0$. From \eqref{eqbd1}, on the event that $\sup_{s\leq T}|c_s^N|\leq K$,
\begin{align*}
&\sqrt{N}\sum_{\epsilon\leq t\leq T,\Delta \Lambda^N_t>0} \tilde V^{\Lambda^N}(t^-,[0,1/N])\\
     &\qquad \leq \tilde C_\epsilon\sum_{i\in \mathbb{N},\epsilon\leq T_i\leq T}\left(\frac{1}{N^{3/2}}+\frac{K}{N}\textrm{erf}\left(\frac{1}{N\sqrt{2(T_{i}-T_{i-1}\vee \epsilon)}}\right)\right)\\
     &\qquad \leq \frac{\tilde C_\epsilon\left( \Lambda_T+\frac{K}{\sqrt{N}}\right)}{\sqrt{N}}+\sum_{i\in\mathbb{N},\epsilon\leq T_i\leq T}\frac{\tilde C_\epsilon K}{N}\textrm{erf}\left(\frac{1}{N\sqrt{2(T_{i}-T_{i-1}\vee \epsilon)}}\right).\end{align*}
The first term converges in probability to zero, as $N\rightarrow\infty.$ By Theorem \ref{bdtheorem}, $\sup_{s\leq T}|c_s^N|$ is stochastically bounded. Further, $K$ is arbitrary. Therefore, to show the result of the current lemma, it is then sufficient to prove that,
$$\left(\sum_{i\in\mathbb{N},\epsilon\leq T_i\leq T}\frac{\tilde C_\epsilon K}{N}\textrm{erf}\left(\frac{1}{N\sqrt{2(T_{i}\vee \epsilon-T_{i-1}\vee \epsilon)}}\right)\right)\mathbbm{1}_{\{\sup_{s\leq T}|c_s^N|\leq K\}}$$
converges in probability to zero.

Since $\frac{N}{2}\partial_x\tilde V^{\Lambda^N}(s^-,0)$ is the jump rate of $\Lambda^N$, for $T_i\geq \epsilon$, \\$\int_{T_{i-1}\vee \epsilon}^{T_{i}}\frac{N}{2}\partial_x\tilde V^{\Lambda^N}(s,0)ds= U_i$ where $U_i$ is a standard exponential random variable, independent of $(\Lambda^N_{s})_{s\leq T_{i-1}}$. By \eqref{eqbd2}, on the event that $\sup_{s\leq T}|c_s^N|\leq K$, it holds that 
$$U_i\leq \tilde C_\epsilon\left(N(T_{i}\vee \epsilon-T_{i-1}\vee \epsilon)+N^{1/2}K\sqrt{T_i\vee \epsilon-T_{i-1}\vee \epsilon}\right).$$
Rearranging this, we determine that for $T_i\geq \epsilon$, 
\begin{align*}\sqrt{T_i-T_{i-1}\vee \epsilon} \geq\frac{ 2U_i}{\tilde C_\epsilon KN^{1/2}+\sqrt{\tilde C_\epsilon^2K^2N+4U_iN}} \geq \frac{U_i\wedge 1}{\tilde C_\epsilon K\sqrt{N}}.\end{align*}

On the event that $\sup_{s\leq T}|c_s^N|\leq K$ ,there are at most $N(\Lambda_T+K/\sqrt{N})$ jumps of $\Lambda^N$ by the time $T$. For all sufficiently large $N$, on this event, we can then bound
\begin{align*}\sum_{i\in\mathbb{N},\epsilon\leq T_i\leq T}\frac{\tilde C_\epsilon K}{N}\textrm{erf}\left(\frac{1}{N\sqrt{2(T_{i}-T_{i-1}\vee \epsilon)}}\right)\leq \sum_{i=0}^{\lceil N(\Lambda_T+1)\rceil}\frac{\tilde C_\epsilon K}{N}\textrm{erf}\left(\frac{\tilde C_\epsilon K}{\sqrt{N}(U_i\wedge 1)}\right).\end{align*}
    The expectation of the right hand side is $\tilde C_\epsilon K(\Lambda_T+1)\mathbb{E}\left(\textrm{erf}\left(\frac{\tilde C_\epsilon K}{\sqrt{N}(U_i\wedge 1)}\right)\right).$ This converges to zero as $N\rightarrow\infty$ by dominated convergence. Therefore, 
    $$\left(\sum_{i\in\mathbb{N},\epsilon\leq T_i\leq T}\frac{\tilde C_\epsilon K}{N}\textrm{erf}\left(\frac{1}{N\sqrt{2(T_{i}\vee \epsilon-T_{i-1}\vee \epsilon)}}\right)\right)\mathbbm{1}_{\{\sup_{s\leq T}|c_s^N|\leq K\}}$$ converges in probability to zero as $N\rightarrow\infty$ As noted earlier, since $K$ is arbitrary, this convergence implies the result of the current lemma.
\end{proof}
\begin{proof}[Proof of Lemma \ref{lto0}]
    We fix arbitrary $\delta>0$, and consider $\epsilon>0$.
    \begin{align*}
    &\limsup_{N\rightarrow\infty}\mathbb{P}(\sqrt{N}\tilde L^N_T>\delta)\\
    &\leq\limsup_{N\rightarrow\infty} \mathbb{P}(\sqrt{N}(\tilde L^N_T-\tilde L_\epsilon^N)>\delta/2)+\limsup_{N\rightarrow\infty}\mathbb{P}(\sqrt{N}\tilde L^N_\epsilon>\delta/2).
    \end{align*}
The first term is exactly zero by Lemma \ref{lto01}. Then sending $\epsilon\downarrow 0$, the second term converges to zero by Lemma \ref{lto02}.
\end{proof}
\section{Functional central limit theorem for the empirical measure}
\label{fffclt}
In this Section we give a simple extension of the functional central limit theorem for $\sqrt{N}(\Lambda^N_t-\Lambda_t)$ to fluctuations of the empirical measure. We define $\Upsilon^N$ by:
$$\Upsilon^N_t(.):=\sqrt{N}\left(\frac{1}{N}\sum_{i=1}^\infty\delta_{X_t^{i,N}}(.)\mathbbm{1}_{\{\tau_i^N>t\}}-\int_.\tilde V^{\Lambda}(t,x)\textrm{d}x\right).$$ 

We write $\mathcal{S}$ to denote the space of Schwartz functions on $\mathbb{R}^+$
%$$\mathcal{S}:=\left\{f\in \mathcal{C}^\infty(\mathbb{R}^+):\forall \alpha,\beta\in\mathbb{N}_0,\sup_{x\geq 0}\left|(1+x)^\alpha \frac{\textrm{d}^\beta f}{dx^{\beta}}(x)\right|<\infty\right\}.$$
We endow this space with the topology generated by the semi-norms $\sup_{x\geq 0}\left|(1+x)^\alpha \frac{\textrm{d}^\beta f}{dx^{\beta}}(x)\right|$, in an identical manner to the standard Schwartz space on $\mathbb{R}$, so that this space is a nuclear, Fr\'echet space. We write $\mathcal{S}'$ for the corresponding dual space of tempered distributions. We then consider $\Upsilon^N$ to be an element of $D([0,T],\mathcal{S}')$ as defined in \cite{Distributionalconvergence}, and endow this space with the topology described therein.

\begin{theorem}
\label{fclthm}
    $\Upsilon^N$ converges in distribution to a process $\Upsilon\in \mathcal{C}([0,T],\mathcal{S}')$. The process $\Upsilon$ is the unique solution to the following SPDE.
    For all $t\geq 0$, and $\phi\in \mathcal{S}$,
    \begin{align*}
        (\Upsilon_t,\phi)=&(\Upsilon_0,\phi)+\frac{1}{2}\int_0^t(\phi'',\Upsilon_s)ds-\int_0^t(\phi',\Upsilon_s)\textrm{d}\Lambda_s\\
        &-c_t(\phi',\tilde V^\Lambda_t)+\int_0^tc_s(\phi',\partial_s\tilde V^\Lambda_s)ds-\phi(0)c_t+S(t,\phi).
    \end{align*}
In this equation, $\Lambda$ is the free boundary in \eqref{alternativetomvr} for the initial data $g$ and $\tilde V^\Lambda$ is the corresponding particle distribution. The process $c$ is the unique solution to $c_t=\delta \Gamma(\Lambda,c)_t+M_t$ for $t\leq T$, where $M$ is a Gaussian process as defined in Theorem \ref{clt}. The processes $\Upsilon_0\in \mathcal{S}'$, $S\in \mathcal{C}([0,T],\mathcal{S}')$, and $M\in \mathcal{C}([0,T])$ are jointly Gaussian with mean 0. Their covariances are given by
%$S$ and $\Upsilon_0$ are independent, while the covariance between other terms can be explicitly computed in terms of $g$ and $\Lambda$.
%Could leave in the covariance computations
%The process $c$, may be defined in terms of $M$ as in Theorem \ref{clt}.  covariances given by:
        \begin{align*}
        &\mathbb{E}(S(t,\phi_1)S(s,\phi_2))=\int_0^{t\wedge s}(\tilde V^\Lambda_r,\phi_1'\phi_2')dr,\\
        &\mathbb{E}(S(t,\phi_1)(\Upsilon_0,\phi_2))=0,\\
        &\mathbb{E}((\Upsilon_0,\phi_1)(\Upsilon_0,\phi_2))=\int_0^\infty  \phi_1(x)\phi_2(x)g(x)\textrm{d}x,\\
        %&\mathbb{E}((\Upsilon_0,\phi)c_t)= \int_0^t\mathbb{E}((\Upsilon_0,\phi)c_s)\mathbb{E}(g(W_s+\Lambda_s)|\tau_t=s)f_{\tau_t}(s)\textrm{d}s\\
        %&\qquad\qquad\qquad\quad +\int_0^\infty \phi(x)\mathbb{P}(\tau^x\leq t)g(x)\textrm{d}x-(g,\phi)\Lambda_t,\\
        &\mathbb{E}((\Upsilon_0,\phi)M_t)=\int_0^\infty \phi(x)\mathbb{P}(\tau^x\leq t)g(x)\textrm{d}x, \\
        %&\mathbb{E}(S(l,\phi)c_t)=\int_0^t\mathbb{E}(S(l,\phi)c_l)\mathbb{E}(g(W_l+\Lambda_s)|\tau_t=l)f_{\tau_t}(l)\textrm{d}l\\%&\qquad\qquad\qquad\quad+\int_0^\infty\mathbb{E}\left(\mathbbm{1}_{\{\tau^x\leq t\}}\int_0^{s\wedge \tau^x}\phi'(x+B_l-\Lambda_l)\textrm{d}B_l\right)g(x)\textrm{d}x,\\
        &\mathbb{E}(S(t,\phi)M_s)=\int_0^\infty\mathbb{E}\left(\mathbbm{1}_{\{\tau^x\leq s\}}\int_0^{t\wedge \tau^x}\phi'(x+B_l-\Lambda_l)\textrm{d}B_l\right)g(x)\textrm{d}x,\\
        &\mathbb{E}(M_tM_s)=\Lambda_{t\wedge s},
    \end{align*}
for arbitrary $ \phi,\phi_1,\phi_2\in \mathcal{S}$, \textrm{ and } $t,s\geq 0.$
Here $\tau^x$ is as in \eqref{alternativetomvr}, and $B$ is a standard Brownian motion.

\end{theorem}
Analogously to the proof of Corollary \ref{4.2}, applying the rescaling given in Section \ref{reformulationsec} we immediately obtain a fluctuation result for long time behaviour of the system \eqref{xidefnequation}.
\begin{corollary}
   Suppose that $a<1$ for the particle system \eqref{xidefnequation}.
   Let $\tilde V^{\Lambda}$ be the similarity solution of the supercooled Stefan problem, corresponding to the free boundary $K_a\sqrt{t}$. The measure valued process 
   $$\Upsilon^t_s(.):=t^{1/4}\left(\frac{1}{\sqrt{t}}\sum_{i=1}^\infty{\delta}_{\left\{\frac{W_{ts}^i-\xi_{st}}{\sqrt{t}}\right\}}(.)\mathbbm{1}_{\{\tau_i>ts\}}-\int_. \tilde V^{\Lambda}(s,x)\textrm{d}x\right),$$
   when viewed as an element of $D([0,T],\mathcal{S'})$
   converges in distribution to a process $\Upsilon\in \mathcal{C}([0,T],\mathcal{S'})$, which solves the SPDE given in Theorem \ref{fclthm}, for the free boundary $\Lambda_t=K_a\sqrt{t}$, and $g\equiv a$.
\end{corollary}

 After applying the previously established bounds on $c^N_t$ the proof of this Theorem \ref{fclthm} follows mainly by standard arguments. We show tightness of the process $\Upsilon^N$, due to the bounds on $c^N_t$, then use Ito's formula to derive the SPDE solved by the limiting process. To show uniqueness, we then determine the joint distribution of $(\Upsilon_0,S,c)$. This is done by considering an auxiliary particle system.
\subsection{Properties of $\Upsilon^N$}

We couple the particle system to an auxiliary system given by $\tilde X^{i,N}= X^{i,N}(\Lambda)$ and $\tilde \tau_i^N=\tau_i^N(\Lambda)$, for ${i\in\mathbb{N}}$, as used in defining $\Gamma_N(\Lambda)$ in Section \ref{deterministic barriers}. We then set $\tilde \Lambda^N_t=\Gamma_N(\Lambda)_t$, $\tilde c_t^N=\sqrt{N}(\tilde \Lambda^N_t-\Lambda_t)$ and define
$$\tilde\Upsilon^N_t(.):=\sqrt{N}\left(\frac{1}{N}\sum_{i=1}^\infty\delta_{\tilde X_t^{i,N}}(.)\mathbbm{1}_{\{\tilde \tau_i^N>t\}}-\int_.\tilde V^{\Lambda}(t,x)\textrm{d}x\right).$$
Since the particles in this auxiliary system are entirely independent, it is simple to control $\tilde \Upsilon^N$ and $\tilde c^N$ in order to determine their joint limit. We shall compare $\tilde \Upsilon^N$ to $\Upsilon^N$ in order to determine properties of $\Upsilon^N$.

The main noise term which will appear in the limit corresponds to the following processes $S^N,\tilde S^N\in D((0,\infty),\mathcal{S}')$ which we define as
$$S^N(t,\phi):=\frac{1}{\sqrt{N}}\sum_{i=1}^\infty \int_0^t\phi'(X_s^{i,N})\mathbbm{1}_{\{\tau_i^N>s\}}\textrm{d}B_s^{i,N},$$
$$\tilde S^N(t,\phi):=\frac{1}{\sqrt{N}}\sum_{i=1}^\infty\int_0^t \phi'(\tilde X_s^{i,N})\mathbbm{1}_{\{\tilde \tau_i^N>s\}}\textrm{d}B_s^{i,N}.$$

%We extend $\Upsilon^N_0$ as an element of $D([0,T],\mathcal{S}')$ by setting $(\Upsilon_0^N)_s=\Upsilon_0^N$ for all $s\leq T$. We then view 
%$(S^N,\tilde S^N,\Upsilon^N_0,\tilde c^N)$ as an element of $D([0,T],\mathcal{S}'\times \mathcal{S}'\times \mathcal{S}'\times \mathbb{R})$. We endow the space $E'=\mathcal{S}'\times \mathcal{S}'\times \mathcal{S}'\times \mathbb{R}$ with the product topology, and endow the space $D([0,T],E')$ with the corresponding topology on as given in \cite{Distributionalconvergence}.

The following results are simple consequences of the definitions of the processes, and the convergence of $c^N$. The proofs are included in Appendix \ref{fffcltappen}.
\begin{lemma}
\label{jointdist}
    $(S^N,\tilde S^N,\Upsilon^N_0,\tilde c^N)$ when viewed with the product topology converges in distribution to $(S,S,\Upsilon_0,\tilde c)$. The triple $(S,\Upsilon_0,\tilde c)$ is Gaussian, and $(S,\Upsilon_0,\tilde c)$ has the covariance structure of $(S,\Upsilon_0,M)$ as given in Theorem \ref{fclthm}.
\end{lemma}

\begin{lemma}
\label{tightnessUpsilon^N}
    For $\phi \in \mathcal{S}$, and $T\geq 0$, $(\Upsilon_.^N,\phi)$ is tight in $D([0,T],\mathbb{R})$, endowed with the uniform topology. Any limit is in $\mathcal{C}([0,T],\mathbb{R})$.
\end{lemma}
\begin{corollary}
\label{tightUpsilon}
    $\Upsilon_.^N$ is tight when considered as a process in $D([0,T],\mathcal{S}')$, and any subsequential limit is almost surely an element of $\mathcal{C}([0,T],\mathcal{S}')$.
\end{corollary}

\begin{proof}
    Tightness follows from \cite[Theorem 3.1]{Distributionalconvergence} and the result of Lemma \ref{tightnessUpsilon^N}. Further, since the uniform topology is used in Lemma \ref{tightnessUpsilon^N}, we also observe that for any subsequential limit $\Upsilon$, almost surely, for all $\lambda \in \mathbb{Q}^+$, $(\Upsilon,e^{-\lambda x})$ must be continuous.
    Suppose that $\Upsilon$ is discontinuous with positive probability, so that a random time $s$ exists, at which the c\'adl\'ag process $\Upsilon$ jumps. Almost surely, at the time $s$ it holds that $(\Delta \Upsilon_s,e^{-\lambda x})=0$ for all $\lambda \in \mathbb{Q}^+$. Therefore $\Delta \Upsilon_s=0$, contradicting the assumption that $s$ was a jump time. 
\end{proof}

The form of the limiting SPDE then appears from an application of Ito's formula and is given in Appendix \ref{fffcltappen}.
\begin{lemma}
\label{itoforspde}
    For any subsequential limit of $(\Upsilon^N,S^N,c^N)$, almost surely this limit $(\Upsilon,S,c)$ satisfies the following equation for all $t\geq 0$, and $\phi\in \mathcal{S}$,
    \begin{align*}
        (\Upsilon_t,\phi)=&(\Upsilon_0,\phi)+\frac{1}{2}\int_0^t(\phi'',\Upsilon_s)ds-\int_0^t(\phi',\Upsilon_s)\textrm{d}\Lambda_s\\
        &-c_t(\phi',\tilde V^\Lambda_t)+\int_0^tc_s(\phi',\partial_s\tilde V^\Lambda_s)ds-\phi(0)c_t+S(t,\phi).
    \end{align*}
\end{lemma}
\begin{remark}
\label{uniquenessremark!}
For a given $( S,c,\Upsilon_0)$ it is simple to check that the solution to the above equation with $\Upsilon_t \in \mathcal{C}([0,T],\mathcal{S'})$ is unique. This may be done for example by taking two solutions $\Upsilon, \hat \Upsilon$ and noting that $P_t(\lambda)=(\Upsilon_t-\hat \Upsilon_t,e^{-\lambda x})$ must solve
$$P_t(\lambda)=\frac{\lambda^2}{2}\int_0^tP_s(\lambda)ds+\lambda \int_0^tP_s(\lambda)\textrm{d}\Lambda_s.$$
This has unique continuous solution $P_t(\lambda)=0$ for all $t\geq 0$.
\end{remark}

\subsection{Joint distribution of limiting terms}
To finish the proof of Theorem \ref{fclthm}, it is sufficient to determine the joint distribution of any subsequential limit $(S,c,\Upsilon_0)$. 
We characterise this distribution by giving $c$ in terms of $\tilde c $, since the distribution of $(S,\Upsilon_0,\tilde c)$ is previously determined in Lemma \ref{jointdist} .

We determine $c$ in terms of $\tilde c$ by comparison of $\tilde c_t^N$ with $M^N_t+\bar M^N_t$. We first need to control over terms of the form $(\Delta \Lambda^N_s)^2$. We define 
$$L_t^N:=N\sum_{s\leq t}(\Delta \Lambda^N_s)^2.$$

\begin{lemma}
\label{Quadratictight}
    For all positive times $t$,
    $$\mathbb{E}(L_t^N)\rightarrow \Lambda_t \textrm { as }N\rightarrow \infty.$$
\end{lemma}
We then determine that $\tilde c_t$ will be equal to $M_t.$
\begin{lemma}
\label{lemmahalfwaytodist}
    For any fixed $t>0$, $$\limsup_{N\rightarrow\infty}N\mathbb{E}\left(\left((\tilde \Lambda^N_t-\Lambda_t)-(M^N_t+\bar M^N_t)\right)^2\right)=0.$$
\end{lemma}

The proofs of these results are lengthy, and are similar to the proof of Theorem \ref{fclt}. We give the proofs in Appendix \ref{fffcltappen}.

\begin{lemma}
\label{lemmaondist}
    For any subsequential limit of $(c^N,\tilde c^N)$, almost surely $c_t=\delta \Gamma(\Lambda,c)_t+\tilde c_t$ for all $t\leq T$.
\end{lemma}
\begin{proof}[Proof of Lemma \ref{lemmaondist}]
We consider $(\tilde c^N,c^N,M^N+\bar M^N,(\tilde c^N-M^N-\bar M^N)_{s\in\mathbb{Q}\cap[0,T]})$, viewed as an element of $D([0,T],\mathbb{R})^{3}\times \mathbb{R}^\infty$, and endowed with the product uniform topology. Taking subsequences and applying the Skorokhod Representation Theorem, we can choose a representation which converges almost surely to $(\tilde c,c,M,(H_s)_{s\in\mathbb{Q}\cap[0,T]})$. From Lemma \ref{lemmahalfwaytodist}, $(H_s)_{s\in \mathbb{Q}}=(0)_{s\in\mathbb{Q}}$. By definition, it holds also that $0=H_s=\tilde c_s-M_s$ for $s\in\mathbb{Q}\cap[0,T]$. Since the limit processes are c\'adl\'ag almost surely, $0=H_s=\tilde c_s-M_s$ for all $s\leq T$ almost surely. Therefore, almost surely, for all $t\leq T$, $\tilde c_t=M_t$. From the arguments of the proof of Theorem \ref{clt}, it follows also that $c_t=\delta \Gamma(\Lambda,c)_t+M_t=\delta\Gamma(\Lambda,c)_t+\tilde c_t$ for all $t$.
\end{proof}
Relabelling the limit $\tilde c$ as $M$ in keeping with the notation of Theorem \ref{fclthm}, we then establish that the joint distribution of any limiting $(\Upsilon_0,S,M,c)$ is exactly as given in Theorem \ref{fclthm}.
\begin{corollary}
\label{jointdist2}
    Consider any subsequential limit $(S,\Upsilon_0,c,M)$ of $(S^N,\Upsilon^N_0,c^N,\tilde c^N)$. Then $(S,\Upsilon_0,M)$ is jointly Gaussian, with covariance as given in Theorem \ref{fclthm}. The process $c$ is the unique solution to $c_t=\delta \Gamma(\Lambda,c)_t+M_t$ for all $t\leq T$.
\end{corollary}
\begin{proof}
    We may consider a further subsequence on which $(S^N,\tilde S^N,\Upsilon^N_0,c^N,\tilde c^N)$ converges to $(S,S,\Upsilon_0,c,M)$. Recall from Lemma \ref{jointdist} that $(S,\Upsilon_0,M)$ are jointly Gaussian, with covariance structure as given in Theorem \ref{fclthm}. By Lemma \ref{lemmaondist}, almost surely $c_t=\delta \Gamma(\Lambda,c)_t+M_t$. 
\end{proof}
From this, Theorem \ref{fclthm} immediately follows.
\begin{proof}[Proof of Theorem \ref{fclthm}]
    From Corollary \ref{tightUpsilon}, and Lemmas \ref{jointdist}, \ref{tightnessUpsilon^N} any subsequence of $(\Upsilon^N,\Upsilon^N_0,S^N,c^N,\tilde c^N)$ viewed with the product topology possesses a convergent subsequence. From Lemma \ref{itoforspde} it holds that any subsequential limit is the unique solution to the SPDE driven by $(\Upsilon_0,S,c,M)$. From Corollary \ref{jointdist2}, we identify that $(S,\Upsilon_0,M)$ is a Gaussian process with the given covariance structure, and $c$ is the unique solution to $c_t=\delta \Gamma(\Lambda,c)_t+M_t$ for $t\leq T$. Therefore, by Remark \ref{uniquenessremark!}, any subsequential limit has the same distribution, and we obtain the convergence in distribution of the entire sequence.
\end{proof}

\section*{Funding}

This publication is based on work supported by the EPSRC Centre for Doctoral
Training in Mathematics of Random Systems: Analysis, Modelling and Simulation (EP/S023925/1). D.G.M. Flynn is funded by the Charles Coulson Scholarship. 
For the purpose of open access, the authors have applied a CC BY public
copyright licence to any author accepted manuscript arising from this
submission.

\bibliography{bible}
\bibliographystyle{siam}

\appendix
\section{Existence of classical solutions}
\label{classical solutions}
\begin{proposition}
    \label{smoothness}
    Suppose that $g$ satisfies the weak feedback condition \hyperlink{W}{(W)}. Then there exists a unique classical solution to the supercooled Stefan problem on the time interval $(0,\infty)$ . This corresponds to the solution of \eqref{alternativetomvr}.
\end{proposition}
\begin{proof}
For $n\in\mathbb{N}$, we use $g^n$ to denote the function $g^n(x)=g(x)\mathbbm{1}_{[1/n,n]}(x)$. By \cite[Theorem 1.1]{PDEresult}, for each $n\in\mathbb{N}$, the supercooled Stefan problem with initial data $g^n=g\mathbbm{1}_{\{[1/n,n]\}}$ possesses a unique classical solution on $(0,\infty)\times \mathbb{R}^+$. By uniqueness, this is also the solution to \eqref{alternativetomvr}. We denote the corresponding free boundary by $\Lambda^n$. This is increasing in $n$.

Arguing as in \cite[Section 4]{https://doi.org/10.1002/mma.1670050107}, we can determine the existence of a constant $\tilde C$ such that $|\dot\Lambda^n_t|\leq \tilde Ct^{-1/2}$ for all sufficiently large $n$. Passing to the limit yields a solution to \eqref{alternativetomvr} such that $\Lambda$ is Lipschitz on $[\delta,\infty)$ for any positive $\delta$. This must be the unique solution to \eqref{alternativetomvr}. Since $\Lambda$ is Lipschitz on $[\delta,\infty)$, it is standard that $\tilde V^{\Lambda}(t,x-\Lambda_t)$ yields a classical solution to the Stefan problem on $(0,\infty)\times \mathbb{R^+}$. This can be established by working on $[0,T]\times [0,M]$ for arbitrary $T>0$, $M>\Lambda_T$, applying the arguments of \cite[Section 4]{https://doi.org/10.1002/mma.1670050107}, using the Lipschitz bound on $\Lambda$, and applying standard heat kernel bounds on the effect that the data $g\mathbbm{1}_{[n,\infty)}$ has upon the behaviour of the solution density in $[0,M]$.
\end{proof}
\section{Proofs deferred from Section \ref{fclt}}
\label{fcltappen}

\begin{proof}[Proof of Lemma \ref{mbound}]
    The fact that the process considered is a martingale with the given quadratic variation is clear, and the proof is omitted.
    
    To prove the final statement, first recall from Theorem \ref{bdtheorem} that the process $\Lambda^N$ converges to $\Lambda$ uniformly in probability on compact sets. By the Skorokhod Representation Theorem, let us consider versions of $\Lambda^N$ which converge almost surely, uniformly on $[0,T]$, to the process $\Lambda$. Abusing notation, we use the same notation for this representation.
    For arbitrary $\epsilon,\delta>0$, and $t\leq T$ we may bound
\begin{align*}
    &\limsup_{N\rightarrow\infty}\int_0^t\int_0^{\frac{N}{2}\partial_x\tilde V^{\Lambda^N}(s^-,0)}\tilde V^{\Lambda^N}(s^-,[0,\Delta\Lambda_s^N])Q(\textrm{d}w,\textrm{d}s),\\
    &\leq\limsup_{N\rightarrow\infty}\left( \sup_{s\leq \delta}\sup_{x\leq \Delta\Lambda_s^N}\tilde V^{\Lambda^N}(s^-,x)\Lambda^N_\delta+\sup_{s> \delta}\sup_{x\leq \Delta\Lambda_s^N}\tilde V^{\Lambda^N}(s^-,x)\Lambda^N_T\right),\\
    &\leq \limsup_{N\rightarrow\infty}\left(\sup_{s\leq T}\sup_{x>0}\tilde V^{\Lambda^N}(s,x)\Lambda^N_\delta+\sup_{s>\delta}\sup_{x\leq \epsilon}\tilde V^{\Lambda^N}(s,x)\Lambda^N_T\right),\\
    &\leq ||g||_\infty(\Lambda_\delta+\epsilon)+\sup_{s>\delta}\sup_{x\leq 2\epsilon}\tilde V^{\Lambda-\epsilon}(s,x)(\Lambda_T+\epsilon). 
\end{align*}
The second inequality follows since $\sup_{s\leq T}|\Lambda^N_s-\Lambda_s|<\epsilon/2$ for sufficiently large $N$. The final inequality follows since the value of $V^f(t,x)$ is decreasing in $f$, and since $\sup_{s\leq T}|\Lambda^N_s-\Lambda_s|<\epsilon/2$ for all large $N$. 

Recall that $\Lambda$ is the free boundary in a classical solution to the supercooled Stefan problem for times bounded away from zero. Therefore,  as $\epsilon$ decreases to zero, $\sup_{s>\delta}\sup_{x\leq 2\epsilon}\tilde V^{\Lambda-\epsilon}(s,x)$ also decreases to zero. Since $\epsilon$ was arbitrary, we can therefore bound 
$$ \limsup_{N\rightarrow\infty}\int_0^t\int_0^{\frac{N}{2}\partial_x\tilde V^{\Lambda^N}(s^-,0)}\tilde V^{\Lambda^N}(s^-,[0,\Delta\Lambda_s^N])Q(\textrm{d}w,\textrm{d}s)\leq ||g||_\infty\Lambda_\delta.$$
Further, since $\Lambda$ is continuous at $0$, sending $\delta$ down to zero we obtain that 
$$\limsup_{N\rightarrow\infty}\int_0^t\int_0^{\frac{N}{2}\partial_x\tilde V^{\Lambda^N}(s^-,0)}V^{\Lambda^N}(s^-,[0,\Delta\Lambda_s^N])Q(\textrm{d}w,\textrm{d}s)=0.$$
\end{proof}
\begin{proof}[Proof of Corollary \ref{mt}]
    By the result of Theorem \ref{bdtheorem}, $\Lambda^N$ converges uniformly in probability to $\Lambda$ on compact sets.
    Applying the result of Lemma \ref{mbound}, and noting that $\Gamma(\Lambda^N)_t$ converges in probability to $\Gamma(\Lambda)_t=\Lambda_t$, it follows that the process $M_t^N$ satisfies the conditions of \cite[Theorem 7.1.4]{11}. The result then follows.
\end{proof}
\begin{proof}[Proof of Lemma \ref{isamg}]
    Fix arbitrary $t$ and $N$. We write $F(x)=N\tilde V^{\Lambda^N}(t^-,[0,x])$. We rewrite $$\Delta\Lambda_t^N=\inf\{x:P_{N\Lambda^N_{t^-}}(F(x))<Nx-1\}.$$ Therefore $$F(\Delta\Lambda^N_t)=\inf\{x:P_{N\Lambda^N_{t^-}}(x)<NF^{-1}(x)-1\}.$$ We then observe that, conditional on the occurrence of a jump at time $t$, and on $(\Lambda^N_{s})_{s\in[0,t)}$, $F(\Delta\Lambda^N_t)$ is a stopping time under the natural filtration of the Poisson process $P_{N\Lambda^N_{t^-}}
    (x)$. This stopping time is exactly such that $$\Delta\Lambda^N_t-\frac{1}{N}=\frac{1}{N}P_{N\Lambda^N_{t^-}}(F(\Delta\Lambda^N_t)).$$
   
   Since $V^{\Lambda^N}(t^-,x)\leq ||g||_\infty$, and since points are placed via a Poisson point process, $F(\Delta\Lambda_t^N)$ is bounded above by a scaled version of the busy period of an $M/D/\infty$ queue with arrival rate $||g||_\infty$ and service time one. This is a Borel-Tanner random variable. It is then simple to see that the stopping time $F(\Delta \Lambda^N_t)$ has finite expectation. Further, the Poisson process has bounded increments. Therefore, writing $\mathbb{E}_{t^-}$ for the expectation when conditioning as stated, and applying the Optional Stopping Theorem
    \begin{align*}
        \mathbb{E}_{t^-}\left(\Delta\Lambda^N_t-\frac{1}{N}\right)&=\frac{1}{N}\mathbb{E}_{t^-}\left(\int_0^{F(\Delta\Lambda^N_{t})}dP_{N\Lambda^N_{t^-}}(x)\right),\\&=\frac{1}{N}\mathbb{E}_{t^-}(F(\Delta\Lambda^N_t)),\\&=\mathbb{E}_{t^-}(\tilde V^{\Lambda^N}(t^-,[0,\Delta\Lambda^N_{t}])).
    \end{align*}
Since the collection $(P_i)_{i\in\mathbb{N}}$ is independent of the Poisson random measure $Q$, and since $t$ above was arbitrary, we can then conclude that for any $t\geq l\geq 0$
    $$\mathbb{E}\left(\int_l^t\int_0^{\frac{N}{2}\partial_x\tilde V^{\Lambda^N}(s^-,0)}\left(\Delta \Lambda^N_{s}-\frac{1}{N}-V(s^-,[0,\Delta\Lambda^N_s])\right)Q(\textrm{d}w,\textrm{d}s)\middle |(\Lambda^N_{h})_{h\in[0,l]}\right)=0.$$
\end{proof}

\begin{proof}[Proof of Corollary \ref{mg2conv}]
    The quadratic variation of this martingale is given by
    \begin{align*}&N\int_0^t\int_0^{\frac{N}{2}\partial_x\tilde V^{\Lambda^N}(s^-,0)}\left(\Delta\Lambda^N_s-1/N-V^{\Lambda^N}(s^-,[0,\Delta\Lambda^N_s])\right)^2Q(\textrm{d}w,\textrm{d}s)\\
    &\qquad =\frac{1}{N}\int_0^t\int_0^{\frac{N}{2}\partial_x\tilde V^{\Lambda^N}(s^-,0)}\left(\int_0^{NV^{\Lambda^N}(s^-,[0,\Delta\Lambda^N_s])}(P_{N\Lambda^N_{s^-}}(\textrm{d}l)-\textrm{d}l)\right)^2Q(\textrm{d}w,\textrm{d}s).
    \end{align*}

Fix arbitrary $\delta,\epsilon>0$. Since $\Lambda$ is a classical solution to the supercooled Stefan problem on $(0,T]$, we may choose $\delta_\epsilon$ such that $\sup_{\delta<s\leq T}\sup_{x\leq 2\delta_\epsilon}\tilde V^{\Lambda-\delta_\epsilon}(s,x)<\epsilon$. 
We consider only the event that $\sup_{s\leq T}|\Lambda^N_s-\Lambda_s|<\delta_\epsilon$, so that $N\tilde V^{\Lambda^N}(t^-,x)\leq \epsilon$ for $x\leq \Delta\Lambda^N_t$, in the same manner as the proof of Lemma \ref{mbound}.
On this event, for arbitrary $t>\delta$, $N\tilde V^{\Lambda^N}(t^-,[0,\Delta\Lambda^N_t])$ is equal to the stopping time $$\tau_{t^-}^{\epsilon}=\inf\{x:P_{N\Lambda^N_{t^-}}(x)<NF^{-1}(x)-1\}\wedge \inf\{ x:F(x)>N\epsilon x\}.$$ Here we have written $F(x)$ defined as in the proof of Lemma \ref{isamg} above.

As is stated in the proof of Lemma \ref{isamg}, this stopping time has finite expectation and the Optional Stopping Theorem may be applied.
We can therefore bound
\begin{align*}
&\mathbb{E}(([\bar M^N]_t-[\bar M^N]_\delta) \mathbbm{1}_{\{\sup_{s\leq T}|\Lambda^N_s-\Lambda_s|<\delta_\epsilon\}})\\
&\leq \mathbb{E}\left(\int_\delta^t\int_0^{\frac{N}{2}\partial_x\tilde V^{\Lambda^N}(s^-,0)}\mathbb{E}_{s^-}(\tilde V^{\Lambda^N}(s^-,[0,\Delta\Lambda^N_{s}
\wedge \inf\{x:\tilde V^{\Lambda^N}(s,[0,x])>\epsilon x\}]
))Q(\textrm{d}w,\textrm{d}s)\right),\\
&\leq \mathbb{E}\left(\epsilon\int_\delta^t\int_0^{\frac{N}{2}\partial_x\tilde V^{\Lambda^N}(s^-,0)}\Delta\Lambda^N_{s}
Q(\textrm{d}w,\textrm{d}s)\right),\\
&=\mathbb{E}(\epsilon \Lambda^N_t).
\end{align*}
We have written $\mathbb{E}_{s^-}$ for the conditional expectation used in the proof of Lemma \ref{isamg}.
%Rather than this epsilon we could instead try and bound \Delta\Lambda^N_s-1/N.
On the time interval $[0,\delta]$, we can instead bound the expected quadratic variation
\begin{align*}
&\mathbb{E}([\bar M^N]_\delta\mathbbm{1}_{\{\sup_{s\leq T}|\Lambda^N_s-\Lambda_s|<\delta_\epsilon\}})\\
&\leq \mathbb{E}\left(\int_0^\delta\int_0^{\frac{N}{2}\partial_x\tilde V^{\Lambda^N}(s^-,0)}\mathbb{E}_{s^-}(\tilde V^{\Lambda^N}(s^-,[0,\Delta\Lambda^N_{s}
\wedge \inf\{x:\tilde V^{\Lambda^N}(s,[0,x])>\epsilon x\}]
))Q(\textrm{d}w,\textrm{d}s)\right),\\
&\leq \mathbb{E}\left(C\int_0^\delta\int_0^{\frac{N}{2}\partial_x\tilde V^{\Lambda^N}(s^-,0)}\Delta\Lambda^N_{s}
Q(\textrm{d}w,\textrm{d}s)\right),\\
&=\mathbb{E}( C\Lambda^N_\delta).
\end{align*}
Applying Markov's inequality, we obtain that for arbitrary $x>0$,
\begin{align*}
    &\limsup_{N\rightarrow\infty}\mathbb{P}([\bar M^N]_T>x)
    \\
    &\qquad\leq \limsup_{N\rightarrow\infty}\mathbb{P}(\sup_{s\leq T}|\Lambda^N_s-\Lambda_s|>\delta_\epsilon)+\limsup_{N\rightarrow\infty}4\epsilon\mathbb{E}(\Lambda^N_T)/x^2+\limsup_{N\rightarrow\infty}4C\mathbb{E}(\Lambda^N_{\delta})/x^2.
\end{align*}
The first term on the right hand side above is zero. Sending $\epsilon\downarrow 0$, the second term converges to zero. Applying the bound from Theorem \ref{bdtheorem} we can determine that $\mathbb{E}(\Lambda^N_\delta)\leq \Lambda_\delta+\frac{K}{\sqrt{N}}$ for some constant $K$. Therefore the final term above is less than $C\Lambda_\delta$. Sending $\delta\downarrow 0$, and since $\Lambda$ is a continuous function, we obtain finally that $$\limsup_{N\rightarrow\infty}\mathbb{P}([M^N]_T>x)=0.$$
Since $x$ was arbitrary, applying \cite[Theorem 7.1.4]{11} we conclude that $\bar M^N_t$ converges uniformly to zero, on $[0,T]$. 
\end{proof}

The proof of Lemma  \ref{tightlemma2} follows from the following lemma.
    \begin{lemma}
\label{tightN}
    For any $T\geq h>0$:
    $\sqrt{N}(\Lambda^N_t-\Lambda_t)_{t\in [h,T]}$ is uniformly tight.
\end{lemma}
The proof is deferred until after that of Lemma \ref{tightlemma2}. Assuming this result, the proof of Lemma \ref{tightlemma2} is then simple.

\begin{proof}[Proof of Lemma \ref{tightlemma2}]    
    To show tightness, it is equivalent to show that for any subsequence of $\sqrt{N}(\Lambda^N_t-\Lambda_t)_{t\in [0,T]}$, there is a further subsequence which converges in distribution under the uniform topology on $[0,T]$. 
    
    By Lemma \ref{tightN}, for any sequence of integers growing to infinity, there is a further subsequence, $N_j\uparrow \infty$, such that for any $K\in \mathbb{N}$, $\sqrt{N_j}(\Lambda^{N_j}_t-\Lambda_t)_{t\in [1/K,T]}$ converges uniformly in distribution to some continuous process $(c_t(K))_{t\in [1/K,T]}$. By the result of Theorem \ref{bdtheorem}, we may consider a further subsequence on which $\sup_{s\leq 1/m }|c^N_s|$ converges in distribution to some $f(m)$ for each $m\in \mathbb{N}$. 

    We fix a subsequence of integers $N_i\uparrow\infty$ along which the aforementioned convergence holds. Applying the Skorokhod Representation Theorem, we can obtain representations of $\sqrt{N_i}(\Lambda^{N_i}_t-\Lambda_t)_{t\in [0,T]}$, $\sqrt{N_i}(\Lambda^{N_i}_t-\Lambda_t)_{t\in [1/K,T]}$, and $\sup_{s\leq 1/m }|c^{N_i}_s|$ which converge almost surely. From the bounds given in Lemma \ref{lem:aplus} and Lemma \ref{lem:aminus}, we can also observe that $f(m)$ converges in distribution to zero, as $m\rightarrow\infty.$
    We can then take a further representation for which $f(m)$ also converges to zero. The limit processes must satisfy that for any $K\in \mathbb{N}$, $c_t(K)=c_t(K+1)$ on $[1/K,T]$ and that $f(m)\geq \sup_{1/K\leq s\leq 1/m}|c_s(K)|$ for all integers $m$ and $K$, with $m\leq K$. 
     
     We define a process $c$ by $c_t=\liminf_{K\rightarrow \infty}c_t(K)$ for $t>0$ and $c_0=0$. It is simple to check that $\sqrt{N_i}(\Lambda^{N_i}_t-\Lambda_t)$ converges uniformly to $c_t$ almost surely. Fix an arbitrary $\epsilon>0$. We can choose $K$ sufficiently small such that 
     $\sup_{s\leq 1/K}|c_s|<\epsilon/3$ and such that
     $\sup_{s\leq 1/K }|c^N_s|<\epsilon/3$ for all $N_i$ sufficiently large. Further, for all $N_i$ sufficiently large, we may bound $\sup_{t\in [1/K,T]}|\sqrt{N_i}(\Lambda^{N_i}_t-\Lambda_t)-c_t(K)|<\epsilon/3$. Since $\epsilon$ was arbitrary, this proves the convergence.
     
     For any subsequence of integers growing to infinity, we have therefore identified a further subsequence on which $\sqrt{N_i}(\Lambda^{N_i}_t-\Lambda_t)_{t\in [0,T]}$ converges uniformly in distribution. Hence this process is uniformly tight.
     
The second part of Lemma \ref{tightlemma2} follows from the first by rearranging the equation \eqref{rewriteeqtn} to obtain
$$\sqrt{N}(\Gamma(\Lambda^N)_t-\Gamma(\Lambda)_t)=\sqrt{N}(\Lambda^N_t-\Lambda_t)-\bar M^N_t-M^N_t.$$
From the first part of the lemma and Corollaries \ref{mt} and  \ref{mg2conv}
each term on the right hand side of this equation is uniformly tight. Therefore the left hand side of this equation is also uniformly tight.
\end{proof}
%CHange the notation
%Need to remove weak feedback references.
\begin{proof}[Proof of Lemma \ref{tightN}]

For an initial density $g^N$, we shall denote the solution to the particle system \eqref{aligned:Particlealphanocentre} with initial density $g^N$ by $\Lambda^{N,g^N}$. This initial density may be random, as long as it is independent of the Poisson random measure and the collection of Poisson processes used in constructing $\Lambda^{N,g^N}$ given in equation \ref{newrep}. We also denote the minimal solution to \eqref{alternativetomvr}, with initial density $g^N$, by $\Lambda^{g^N}$. We denote $\sqrt{N}(\Lambda^N_t-\Lambda_t)$ by $c^N_t$, and $\sup_{s\leq T}|c_s^N|/\sqrt{N}$ by $m_N$. Finally, in this proof we denote the function $\sup_{y\in(x,x+2m_N)}\tilde V^{\Lambda-m_N}(s,y)$ as $Z^N(s,x)$.

We construct an upper bound on the increments of $c_t^N$. A lower bound can be created in a similar manner. We first bound
\begin{align}
        \Lambda^N_t-\Lambda^N_s-(\Lambda_t-\Lambda_s)
        &\overset{d}{=}\Lambda^{N,\tilde V^{\Lambda^N}(s,.)}_{t-s}-\Lambda^{\tilde V^\Lambda(s,.)}_{t-s},\nonumber\\
        &\leq \Lambda^{N,Z^N(s,.)}_{t-s}-\Lambda^{\tilde V^{\Lambda}(s.)}_{t-s},\nonumber\\
        &\leq (\Lambda^{N,Z^N(s,.)}_{t-s}-\Lambda^{Z^N(s,.)}_{t-s})+(\Lambda^{Z^N(s,.)}_{t-s}-\Lambda^{\tilde V^{\Lambda}(s,.)}_{t-s}).\label{lambdadiff}
    \end{align}
%Since $h>0$, and since $\Lambda$ is a classical solution to the supercooled Stefan problem, there exists $\kappa>0$, $\delta_h,\epsilon_h$ such that if $c\leq \epsilon_h$, $$\mathbb{E}^W\left(\sup_{\{x\leq (\Lambda^{Z^N(s,.)}_{\delta_h}+2\epsilon_h+\sup_{s\leq \delta_h}W_s)\}} Z^N(s,x)\middle |(c_{l}^N)_{l\leq s}\right)\leq 1-\kappa<1.$$
%$W$ is a Brownian motion, which is independent to $\Lambda^N$. Through similar arguments to \cite[Proposition 2.2]{2}, we can then see that there exists some constant $H_h$ such that $\sup_{h\leq s\leq T}\sup_{l\leq \delta_h}\Lambda_l^{Z^N(s,.)}\leq H_h\sqrt{l}.$ Further, through similar arguments to Theorem \ref{bdtheorem}, and inspecting Lemma \ref{lem:aplus}, we can determine that for $\delta\leq \delta_h$, and $x\leq \epsilon_h$,
It is clear that $Z^N(s,.)\leq ||g||_\infty=1-\kappa$ for some $\kappa>0$. Thus, by the same argument as in \cite[Proposition 2.2]{2},
\begin{align*}
\mathbb{P}(\sup_{l\leq \delta}\sqrt{N}|\Lambda^{N,Z^N(s,.)}_{l}-\Lambda^{Z^N(s,.)}_l|>x)&\leq e^{-\frac{\kappa^2 x^2}{2C\Lambda^{Z^N(s,.)}_\delta+2(C+1)x/\sqrt{N}}} 
.\end{align*}
By taking $x=(\Lambda_{\delta}^{Z^N(s,.)})^{1/3}$ in this expression, we observe that
\begin{align}\limsup_{N\rightarrow\infty}\mathbb{P}(\sup_{l\leq \delta}\sqrt{N}|\Lambda^{N,Z^N(s,.)}_{l}-\Lambda^{Z^N(s,.)}_l|>(\Lambda^{Z^N(s,.)}_\delta)^{1/3})\leq e^{-\frac{\kappa^2}{C(\Lambda^{Z^N(s,.)}_\delta)^{1/3}}}.\label{lambdadiff1}\end{align}
This gives a bound on the first term in \eqref{lambdadiff}. 

To consider the second term in \eqref{lambdadiff}, we first determine bounds on $Z^N$. Recall that $\Lambda$ is assumed continuously differentiable on $[h/3,T],$ and therefore $\Gamma(\Lambda-m_N)_l$ is uniformly Lipschitz on the interval $[h/2,T]$, with uniform bounds over sufficiently small $m_N$. Using the definition of $\tilde V^{\Lambda-m_N}$, we can explicitly write 
\begin{align*}
    \tilde V^{\Lambda-m_N}(s,x)=\int_0^\infty\frac{e^{-\frac{(x+\Lambda_s-\Lambda_h-y)^2}{2(s-h/2)}}}{\sqrt{2\pi(s-h/2)}} \tilde V^{\Lambda-m_N}(h/2,y)\textrm{d}y+\int_{h/2}^{s}\frac{e^{-\frac{(x+\Lambda_s-\Lambda_l)^2}{2(s-l)}}}{\sqrt{s-l}}\partial_l\Gamma(\Lambda-m_N)_{l}\textrm{d}l.
\end{align*}

Since $\tilde V^{\Lambda-m_N}\leq ||g||_\infty$, it is standard to check that the derivative in $x$ exists, and is uniformly bounded over $h\leq s\leq T$, so that $\tilde V^{\Lambda-m_N}(s,x)$ is uniformly Lipschitz in $x$ over $h\leq s\leq T$. Therefore, by the definition of $Z^N$ to determine the existence of a constant $D_h$ such that \begin{align}\sup_{h\leq s\leq T}\sup_{x\geq 0}|Z^N(s,x)-\tilde V^\Lambda(s,x)|\leq D_hm_N\label{eqtnv1}\end{align} it is sufficient to show the existence of a constant $C(h,\Lambda)$ such that \begin{align}\sup_{h\leq s\leq T}\sup_{x\geq 0}|\tilde V^{\Lambda-m_N}(s,x+m_N)-\tilde V^{\Lambda}(s,x)|\leq C(h,\Lambda)m_N.\label{eqtnv12}\end{align}

To show this, we first rewrite $\tilde V^{\Lambda-m_N}(s,x+m_N)-\tilde V^{\Lambda}(s,x)$ as the Lebesgue density of the measure
\begin{align*}
    \int_0^\infty\mathbb{P}(y+B_s-\Lambda_s\in \textrm{d}x,0\geq \inf_{l\leq s}( y+B_l-\Lambda_l)>-m_N)g(y)\textrm{d}y
\end{align*}
Since $\Lambda$ is the free boundary in a classical solution to the supercooled Stefan problem on $(0,\infty)$, we can compute that the above density is equal to $k(s,x,m_N)$ where $k$ is defined by
$$k(s,x,m):=\int_0^s\left(\frac{1}{\sqrt{2\pi(s-l)}}\left(e^{-\frac{(x+\Lambda_s-\Lambda_l)^2}{2(s-l)}}-e^{-\frac{(x+\Lambda_s-\Lambda_l+2m)^2}{2(s-l)}}\right)\right)\dot \Lambda_l\textrm{d}l\mathbbm{1}_{\{x>-m\}}.$$
    For a fixed $s> 0$, and $x>-m$, we can differentiate the above in $m$. We obtain that
\begin{align*}\partial_m k(s,x,m)=&\int_0^s\left(\frac{(x+\Lambda_s-\Lambda_l+2m)}{(s-l)\sqrt{2\pi(s-l)}}\left(e^{-\frac{(x+\Lambda_s-\Lambda_l+2m)^2}{2(s-l)}}\right)\right)\dot \Lambda_l\textrm{d}l.
\end{align*}
By bounding $\sup_{h/2\leq l\leq T}\dot\Lambda_l$, and by bounding the bracketed terms uniformly over $l\leq h/2$, we can observe that a constant $C(h,\Lambda)$ exists such that
\begin{align*}\partial_m k(s,x,m)&\leq C(h,\Lambda)\left(1+\int_{h/2}^s\left(\frac{(x+C(h,\Lambda)(s-l)+2m)}{(s-l)\sqrt{2\pi(s-l)}}\left(e^{-\frac{(x+\Lambda_s-\Lambda_l+2m)^2}{2(s-l)}}\right)\right)\textrm{d}l\right),
\\&\leq C(h,\Lambda)\left(1+\sqrt{s-h/2}+\int_{h/2}^s\frac{(x+2m)}{(s-l)^{3/2}}e^{-\frac{(x+2m)^2}{2(s-l)}}\textrm{d}l\right),
\end{align*}
where we have increased the constant $C(h,\Lambda)$ in the second inequality. The right hand side of the final inequality is uniformly bounded over $x+2m\geq 0$. At $m=0$, the value of $k(s,x,0)$ is zero for all $x$. Therefore $k(s,x,m_N)\leq C(h,\Lambda)m_N.$ This establishes \eqref{eqtnv12}, which in turn implies the inequality of equation \eqref{eqtnv1}. Thus, there exists some $D_h$ dependent on $h$ and $\Lambda$ such that $Z^N(s,x)-\tilde V^{\Lambda}(s,x)<D_{h}m_N$ uniformly.

With this upper bound established, for $l\leq h$ and $s\geq h$, we then bound
\begin{align*}
    \Lambda^{Z^N(s,.)}_l-\Lambda_l^{\tilde V^{\Lambda}(s,.)}&\leq \sup_{y\leq l}(\Lambda^{Z^N(s,.)}_y-\Lambda_y^{\tilde V^{\Lambda}(s,.)})\mathbb{E}\left(\sup_{\{x\leq \sup_{y\leq l}B_y+\Lambda_l^{Z^N(s,.)}+2m_N\}}Z^N(s,x)\right)\\&\quad+\mathbb{E}\left(\int_0^{\sup_{y\leq l}\Lambda^{\tilde V^{\Lambda}(s,.)}_l+B_y}
    (Z^N(s,x)-\tilde V^{\Lambda}(s,x))\textrm{d}x\right),\nonumber\\
    &\leq (1-\kappa)\sup_{y\leq l}(\Lambda^{Z^N(s,.)}_y-\Lambda_y^{\tilde V^{\Lambda}(s,.)})+D_hm_N\left(\sqrt{\frac{2l}{\pi}}+\sup_{h\leq s\leq T}\dot \Lambda_s l\right).
\end{align*}
It is simple to rearrange this result to obtain that there exists a constant $D'_h$ such that for $l\leq h$, and $h\leq s\leq T$:
\begin{align}\label{lambdadiff2}\sup_{y\leq l}(\Lambda^{Z^N(s,.)}_y-\Lambda_y^{\tilde V^{\Lambda}(s,.)})\leq D_h'm_N\sqrt{l}.\end{align}

By \eqref{lambdadiff}, \eqref{lambdadiff1}, and \eqref{lambdadiff2}, there exists a constant $C'_h$ such that for any $h>t>0:$
\begin{align*}\sup_{h\leq s\leq T}\limsup_{N\rightarrow\infty }\mathbb{P}(\sup_{l\leq t \wedge (T-s)}(c^N_{l+s}-c_s^N)_+\geq D_h'\sqrt{N}m_N\sqrt{t}+(H_h\sqrt{t})^{1/3})\leq e^{-C'_h/t^{1/6}}.
\end{align*}
An identical argument can be used to bound any decrease in $c^N_t$. Therefore,
\begin{align*}\sup_{h\leq s\leq T}\limsup_{N\rightarrow\infty }\mathbb{P}(\sup_{l\leq t\wedge(T-s) }|c^N_{l+s}-c_s^N|\geq D_h'\sqrt{N}m_N\sqrt{t}+(H_h\sqrt{t})^{1/3})\leq 2e^{-C'_h/t^{1/6}}
.\end{align*}
We then note that if $\sup_{|x-y|\leq \delta, \hspace{0.03cm}x,y\leq T }|c^N_x-c^N_y|>\epsilon$, then $\sup_{l\leq \delta }|c^N_{l+s}-c^N_s|>\epsilon/2$ for some $s\in \delta \mathbb{N}\cap [0,T]$.
Therefore, for any $K,$ $\delta>0$ such that $D_h'K\delta^{1/2}+(H_h\sqrt{\delta})^{1/3}<\epsilon/2$:
\begin{align*}\limsup_{N\rightarrow\infty }\mathbb{P}( \sup_{|x-y|\leq \delta, \hspace{0.03cm}x,y\leq T }|c^N_x-c^N_y|>\epsilon)&\leq \left(\frac{1}{\delta}+1\right)e^{-C_h'\delta^{-1/6}}+\limsup_{N\rightarrow\infty}
\mathbb{P}(\sup_{l\leq T}|c_l^N|\geq K),\\
&\leq \left(\frac{1}{\delta}+1\right)e^{-C_h'\delta^{-1/6}}+2e^{-C_T K^2}.
\end{align*}
By sending $\delta\downarrow0$, this implies that for any $K>0$,
\begin{align*}
& \limsup_{\delta\rightarrow 0}\limsup_{N\rightarrow\infty }\mathbb{P}\left( \sup_{x,y\leq T: |x-y|\leq \delta}|c^N_x-c^N_y|>\epsilon\right)\leq 2e^{-C_TK^2}.
\end{align*}
We then send $K\uparrow\infty$ to obtain
\begin{align*}
& \limsup_{\delta\rightarrow 0}\limsup_{N\rightarrow\infty }\mathbb{P}( \sup_{|x-y|\leq \delta, \hspace{0.03cm}x,y\leq T }|c^N_x-c^N_y|>\epsilon)=0.
\end{align*}
Since $\epsilon>0$ was arbitrary, and from the bound of Theorem \ref{bdtheorem}, \cite[Theorem 7.3]{9} applies to $(c^N_t)_{t\in[h,T]}$, from which the result of the current lemma follows. 
\end{proof}

\begin{proof}[Proof of Lemma \ref{5.7}]
   For any subsequential limit over $N_i\rightarrow\infty$, the Skorokhod Representation Theorem allows us to choose representations of $(c^{N_i}_t,y^{N_i}_t)$ which converge uniformly almost surely. It is then sufficient to check pointwise convergence of $y_t^N$ to $\delta \Gamma(\Lambda,c)_t$. Let $t$ be an arbitrary element of the set $[0,T]$.
    \begin{align*}
    \lim_{i\rightarrow\infty}y_t^{N_i}&=\lim_{i\rightarrow\infty}\sqrt{N_i}(\Gamma(\Lambda^{N_i})_t-\Gamma(\Lambda)_t)\\
    &=\lim_{i\rightarrow\infty}\mathbb{E}^W\left(\sqrt{N_i}\int_{\sup_{s\leq t}{(W_s+\Lambda_s)}}^{\sup_{s\leq t}(W_s+\Lambda_s+\frac{c_s^{N_i}}{\sqrt{N_i}})}g(x)\textrm{d}x\right),\\
    &=\mathbb{E}^W\left(\lim_{i\rightarrow\infty}\sqrt{N_i}\int_{\sup_{s\leq t}{(W_s+\Lambda_s)}}^{\sup_{s\leq t}(W_s+\Lambda_s+\frac{c_s^{N_i}}{\sqrt{N_i}})}g(x)\textrm{d}x\right) %\textrm{ by the DCT}, \\
    %&=\mathbb{E}^W\left(c_{\tau_t} %g\left((W+\Lambda)_{\tau_t}\right)\right), 
\end{align*}
by the Dominated Convergence Theorem. 

We claim that since $\Lambda_s$ is increasing, and differentiable on $(0,t]$, almost surely, $\sup_{s\leq t}(W_s+\Lambda_s)$ is a Lebesgue point of $g$. By Girsanov's Theorem, it is clear that for any $\delta>0$, $\sup_{\delta\leq s\leq t}(W_s+\Lambda_s)$ is absolutely continuous with respect to the Lebesgue measure. Therefore, on the event that  the supremum on $[0,t]$ is obtained at time greater than $\delta$, $\sup_{s\leq t}(W_s+\Lambda_s)$ is a Lebesgue point of $g$ almost surely. Since $\Lambda$ is non-decreasing and continuous the probability that the supremum is attained on $[0,\delta]$ tends to zero as $\delta\downarrow 0$. Therefore $\sup_{s\leq t}(W_s+\Lambda_s)$ is a Lebesgue point of $g$ almost surely.
%Is this obvious. The supremum on $\delta<s<t$ is by Girsanov. The probability that we aren't in this set can be made arbitrarily small by sending delta down to zero.
Since $c^N$ converges uniformly to a continuous process $c$, the supremum of $\sup_{s\leq t}(W_s+\Lambda_s+\frac{c_s^N}{\sqrt{N}})$ must occur arbitrarily close to time $\tau_t$, as $N\rightarrow\infty$.
By the Lebesgue differentiation theorem, it then holds that almost surely under the measure $\mathbb{E}^W$,
$$\lim_{i\rightarrow\infty}\sqrt{N_i}\int_{\sup_{s\leq t}{(W_s+\Lambda_s)}}^{\sup_{s\leq t}(W_s+\Lambda_s+\frac{c_s^{N_i}}{\sqrt{N_i}})}g(x)\textrm{d}x=g\left(\sup_{s\leq t}(W_s+\Lambda_s)\right)c_{\tau_t},$$
where $\tau_t$ is the time that $W_s+\Lambda_s$ reaches its maximum in $[0,t]$. Therefore 
 $$\lim_{i\rightarrow\infty}y_t^{N_i}=\mathbb{E}^W\left(g\left(\sup_{s\leq t}(W_s+\Lambda_s)\right)c_{\tau_t}\right).$$
The computations to see that this is identical to $\delta\Gamma(\Lambda,c)_t$ are exactly the same as the calculations for $y^{N_i}t$, upon replacing $c^{N_i}$ in the equation for $y_t^{N_i}$ with the limit $c$.

%I think further explanation is unnecessary.
%In particular, for any $\epsilon>0$, and sufficiently large (random) $N$, $$\sup_{s\leq t}(W_s+\Lambda_s+\frac{1}{\sqrt{N}}(c_s^N))\leq W_{\tau_t}+\Lambda_{\tau_t}+\sup_{s\in (\tau_t-\epsilon,\tau_t+\epsilon)}\frac{c_s^N}{\sqrt{N}}$$

%Taking the limit, the term in the expectation, can be bounded above by $c_s\times \partial^+_x\int_0^xg(y)\textrm{d}y|_{x=\Lambda_{\tau_t}+W_{\tau_t}}$. Since $\tau_t$ has a density, almost surely $\Lambda_{\tau_t}+W_{\tau_t}$ is at a point where the derivative $\partial_x\int_0^xg(y)\textrm{d}y$ exists.

%To lower bound, we instead have, 
%$$\sup_{s\leq t}(W_s+\Lambda_s+\frac{1}{\sqrt{N}}(M_s^N+c_s^N))\geq W_{\tau_t}+\Lambda_{\tau_t}+\frac{1}{\sqrt{N}}(M_{\tau_t}^N+c_{\tau_t}^N),$$
%and we may then argue as for the upper bound.

\end{proof}
\section{Proofs deferred from Section \ref{fffclt}}
\label{fffcltappen}
\begin{proof}[Proof of Lemma \ref{jointdist}]
    The convergence of $(\tilde S^N,\Upsilon^N_0,\tilde c^N)$ to a Gaussian process follows by rewriting these processes as sums of $N$ i.i.d. terms and applying the Central Limit Theorem. By Poisson thinning, and with slight abuse of notation, we relabel the set $(W_t^{i,N})_{i\in\mathbb{N}}$ as $(W_t^{i,k})_{i\in\mathbb{N},k=1,...,N}$. These satisfy that $W_t^{i,k}=W_0^{i,k}+B_t^{i,k}$ where $(W_0^{i,k})_{i\in\mathbb{N}}$ are the points of $N$ independent Poisson processes of intensity $g$, and $(B^{i,k})_{i\in\mathbb{N}}$ are $N$ independent collections of independent Brownian motions. We can then rewrite
    \begin{align*}
        \tilde S^N(t,\phi) &=\frac{1}{\sqrt{N}}\sum_{k=1}^N\left(\sum_{i=1}^\infty\int_0^t\phi'(W_s^{i,k}-\Lambda_s)\mathbbm{1}_{\{\inf_{l\leq s}(W_l^{i,k}-\Lambda_l)>0\}}\textrm{d}B_s^{i,k}\right), \\
        \Upsilon^N_0(\phi) &=\frac{1}{\sqrt{N}}\sum_{k=1}^N\left(\sum_{i=1}^\infty\phi(W_0^{i,k})-\int_0^\infty \phi(x)g(x)\textrm{d}x\right),\\
        \tilde c^N_t&=\frac{1}{\sqrt{N}}\sum_{k=1}^N\left(\sum_{i=1}^\infty\mathbbm{1}_{\{\inf_{s\leq t}(W_s^{i,k}-\Lambda_s)\leq 0\}}-\Lambda_t\right).
    \end{align*}

    From simple adaptations to the arguments in \cite[Theorem 3.1]{Distributionalconvergence} and \cite[Proposition 5.2]{Distributionalconvergence}, in order to prove the convergence in distribution of $(\tilde S^N,\Upsilon^N_0,\tilde c^N)$ to the claimed limit, it is sufficient to show that $\tilde c^N$ is tight, that $\tilde S^N(.,\phi)$ is tight for all $\phi\in \mathcal{S}$, and that for any collections $(\psi_l)_{l=1}^m,(\phi_i)_{i=1}^n\subset \mathcal{S}$, and $(s_i)_{i=1}^n,(t_j)_{j=1}^z\in [0,T]$, 
    $$((\tilde S^N(s_i,\phi_i))_{i=1}^n,(\tilde c^N_{t_j})_{j=1}^z,(\Upsilon_0^N(\psi_l))_{l=1}^m)$$ converges to a Gaussian distribution with the appropriate covariance structure.
    
    The joint distribution can be found by application of the multivariate Central Limit Theorem, together with computations of the covariance terms between $\Upsilon_0^1,$ $\tilde c^1$ and $\tilde S^1$, which are identical to the covariance terms between $\Upsilon_0,$ $\tilde c$ and $S$. We will later compute explicitly $\mathbb{E}(\tilde S^1(t,\phi_1)\tilde S^1(s,\phi_2))$. Other covariance terms can be computed similarly. The tightness of the process $\tilde S^N(.,\phi)$ can be established by showing the existence of a constant $\tilde C$ such that for all $h>0$,
    %\begin{align}\label{tightss}
    %&\sup_{t\leq T}\mathbb{E}((S^1((t+s)\wedge T,\phi)-S^1(t,\phi))^2)\leq \tilde Cs,
    %\end{align}
    %and 
    \begin{align}\label{tights}
    &\sup_{t\leq T}\mathbb{E}((\tilde S^1((t+h)\wedge T,\phi)-\tilde S^1(t,\phi))^4)<\tilde Ch^2.
    \end{align}
Applying Fubini's Theorem, we determine that
    \begin{align*}
    &\mathbb{E}(\left(\tilde S^N((t+h)\wedge T,\phi)-\tilde S^N(t,\phi)\right)^4)\\
    &\qquad\leq 3\mathbb{E}((\tilde S^1((t+h)\wedge T,\phi)-\tilde S^1(t,\phi))^2)^2+\frac{1}{N}\mathbb{E}((\tilde S^1((t+h)\wedge T,\phi)-\tilde S^1(t,\phi))^4),\\
    &\qquad \leq 3\tilde Ch^2+\frac{\tilde Ch^2}{N}.
    \end{align*}
The Kolmogorov tightness criterion can then be applied, since $S^N(0,\phi)=0$.

We now provide a complete proof that for arbitrary $\phi\in\mathcal{S}$ equation \eqref{tights} is satisfied, together with the computation of the covariance terms for the process $S^1$. Since the fourth moment is then finite, we can apply the central limit theorem to any finite collection  $(\tilde S^N(t_i,\phi_i))_{i=1}^n$.
Analogous arguments can be used to show that the central limit theorem can also be applied to $\Upsilon_0(\phi)$, and to $c$. Furthermore, computations similar to those that will be used for the covariance of $\tilde S^1$ can also be applied to $\Upsilon_0^1$ and $\tilde c^1$ to rewrite the covariance terms in the manner of Theorem \ref{fclthm}. The corresponding computations are omitted, since they follow completely analogous arguments to those for $\tilde S^N$.
    
    Let $\phi\in\mathcal{S}$ be arbitrary. To apply the central limit theorem to $\tilde S^1(t,\phi)$, it is sufficient to show that the second moment is bounded. We will instead show that the 4th moment of the increments is suitably bounded, thereby also showing the result of \eqref{tights}. Let $t,h\geq 0$ be arbitrary. We have
    %Change this to 4th moment
    \begin{align*}
         &\mathbb{E}\left(|\tilde S^1((t+h)\wedge T,\phi)-\tilde S^1(t,\phi)|^4\right)\\
         &\qquad \leq \left(\sum_{i=1}^\infty\Bigg(\mathbb{E}\left(\left|\int_t^{(t+h)\wedge T}\phi'(W_s^{i,1}-\Lambda_s)\mathbbm{1}_{\{\inf_{l\leq s}(W_l^{i,1}-\Lambda_l)>0\}}\textrm{d}B_s^{i,1}\right|^4\right)^{1/4}\right)^4
    \end{align*}
where the inequality follows by taking the absolute values in the sum defining $\tilde S^1(t,\phi)$, then applying Fubini's Theorem and the Cauchy-Schwarz inequality twice. Continuing, we have 
\begin{align*}
    &\leq \left(\sum_{i=1}^\infty \mathbb{E}\left(\left(\int_t^{(t+h)\wedge T}|\phi'(W_s^{i,1}-\Lambda_s)|^2\textrm{d}s\right)^2\right)^{1/4}\right)^4,\\
    &\leq h^2\left(\sum_{i=1}^\infty\mathbb{E}\left(\sup\left\{\phi'(x)^4:W_0^{i,1}-\Lambda_T-\sup_{s\leq T}|B_s^{i,1}|\leq x\leq W_0^{i,1}+\sup_{s\leq T}|B_s^{i,1}|\right\}\right)^{1/4}\right)^4,\\
    &\leq \tilde Ch^2\left(\sum_{i=1}^\infty \mathbb{E}\left(1\wedge (W_0^{i,1}-\Lambda_T-\sup_{s\leq T}|B_s^{i,1}|)_+^{-100}\right)^{1/4}\right)^4.
\end{align*}
The first inequality follows from the Burkholder-Davis-Gundy inequality, and the final inequality follows since $\phi\in \mathcal{S}$. We can decrease the value of $W_0^{i,1}$ by considering points placed by a Poisson point process of intensity one, rather than $g$. The random variable $W_0^{i,1}$ is therefore lower bounded by a Gamma$(i,1)$ random variable. By applying standard tail bounds for the Gamma distribution and for the Gaussian, we can then bound the probability that $W_0^{i,1}-\sup_{s\leq T}|B_s^{i,1}|\leq i/4$ by $e^{-\tilde Ci}$ for some constant $\tilde C$. Substituting this into the above equation, it can then be seen that the infinite sum is finite. This is exactly the result of \eqref{tights}.
Finally, we compute the covariance terms. For arbitrary $0\leq s,t\leq T$ and $\phi_1,\phi_2\in \mathcal{S}$, we have
\begin{align*}
    &\mathbb{E}(\tilde S^1(t,\phi_1)\tilde S^1(s,\phi_2))\\
    &\qquad=\sum_{i=1}^\infty\sum_{j=1}^\infty\mathbb{E}\Bigg(\left(\int_0^t\phi'_1(W_l^{i,1}-\Lambda_l)\mathbbm{1}_{\{\inf_{h\leq l}(W_h^{i,1}-\Lambda_h)>0\}}\textrm{d}B_l^{i,1}\right),\\
    &\quad\quad\quad \quad\quad\quad\quad\quad\times \Bigg(\int_0^s\phi'_2(W_l^{j,1}-\Lambda_l)\mathbbm{1}_{\{\inf_{h\leq l}(W_h^{j,1}-\Lambda_h)>0\}}\textrm{d}B_l^{j,1}\Bigg)\Bigg),\\
    &\qquad=\sum_{i=1}^\infty\mathbb{E}\left(\int_0^{t\wedge s}\phi'_1(W_l^{i,1}-\Lambda_l)\phi'_2(W_l^{i,1}-\Lambda_l)\mathbbm{1}_{\{\inf_{h\leq l}(W_h^{i,1}-\Lambda_h)>0\}}\textrm{d}l\right),\\
    &\qquad=\int_0^{t\wedge s}\left(\int_0^\infty\mathbb{E}(\phi'_1(x+B_l-\Lambda_l)\phi'_2(x+B_l-\Lambda_l)\mathbbm{1}_{\{\tau^x\geq l\}})g(x)\textrm{d}x\right)\textrm{d}l,\\
    &\qquad=\int_0^{t\wedge s}(\tilde V^\Lambda_l,\phi'_1\phi'_2)dl.
    \end{align*}
The first equality follows by Fubini's Theorem. The second follows by the independence of the Brownian motions $(B^{i,1})_{i\in\mathbb{N}}$. The third follows by applying Fubini's Theorem twice and using the fact that the points $W^{i,1}_0$ are the points of a Poisson process of intensity $g$, which is independent of the collection $(B^{i,1})_{i\in\mathbb{N}}$. The final equality follows by the definition of $\tilde V^\Lambda$. 

As previously noted, these arguments can then be repeated for $\Upsilon^N$ and $\tilde c^N$, to determine that $(\tilde S^N,\Upsilon^N_0,\tilde c^N)$ converges in distribution to the stated limit. 

To prove the convergence of $(\tilde S^N,S^N,\Upsilon^N_0,\tilde c^N)$ to the desired limit, it is sufficient to check that $S^N-\tilde S^N$ converges in distribution to zero. To do this, it is sufficient to check that for all $\phi\in \mathcal{S}$, $\sup_{s\leq T}|S^N(t,\phi)-\tilde S^N(t,\phi)|$ converges in probability to zero. By \cite[Theorem 3.1]{Distributionalconvergence}, any subsequence of $S^N-\tilde S^N$ would then have a convergent subsequence, and any limit of such a subsequence would be zero by an analogous argument to that of Corollary \ref{tightUpsilon}. We then observe that $S^N(t,\phi)-\tilde S^N(t,\phi)$ is a martingale, whose quadratic variation can be shown to converge in probability to zero, by Theorem \ref{bdtheorem}. Therefore, by \cite[Theorem 7.1.4]{11}, $\sup_{t\leq T}|S^N(t,\phi)-\tilde S^N(t,\phi)|$ converges to zero.

Finally, we need to check the covariances of the limit processes. Since we consider i.i.d. sums of terms with finite second moments, the covariance of the limit processes is exactly the same as the covariance with $N=1$, which can be computed as noted earlier.
\end{proof}

To bound several terms from $\Upsilon^N$, it is necessary to determine the convergence of the particle system under the law of large numbers scaling. This can be immediately deduced from the results of Theorem \ref{bdtheorem} and we state it here for convenience. We define
$$ \mu^N_t:=\frac{1}{N}\sum_{i=1}^\infty\delta_{X_t^{i,N}}(.)\mathbbm{1}_{\{\tau_i^N>t\}}.$$

\begin{lemma}
\label{particlesystemconvergencewithoutscaling}
   $\mu^N$ converges uniformly in probability on $D([0,T],\mathcal{S}')$ to the measure valued process $(\tilde V^{\Lambda}(s,x)\textrm{d}x)_{s\in [0,T]}$.
\end{lemma}

One error term in $\Upsilon^N$, that will be necessary to control is given by  $A^N\in D([0,T],\mathcal{S}')$ defined by
\begin{align*}
A^N_t(\phi):=&\frac{1}{\sqrt{N}}\sum_{i=1}^\infty \sum_{s\leq t}\left(\phi(X_{s^-}^{i,N}-\Delta\Lambda^N_s)-\phi(X_{s^-}^{i,N})+\phi'(X_{s^-}^{i,N})\Delta \Lambda^N_s\right)\mathbbm{1}_{\{\tau_i^N\geq s\}}\\&\quad+\sqrt{N}\phi(0)\Lambda^N_t-\frac{1}{\sqrt{N}}\sum_{i=1}^\infty \sum_{s\leq t}\phi(X_{s^-}^{i,N})\mathbbm{1}_{\{\tau_i^N=s\}}\\
&\quad-\sum_{s\leq t}\Delta \Lambda^N_s\Delta (\mu^N_s,\phi').
\end{align*}
\begin{lemma}
\label{anlemma}
$A^N$ converges in distribution to zero.
\end{lemma}
\begin{proof}[Proof of Lemma \ref{anlemma}]

From \cite[Theorem 3.1, Proposition 5.2]{Distributionalconvergence}, it is sufficient to check that $\sup_{t\leq T}|A^N(\phi)_t|$ converges in probability to zero, for arbitrary $\phi\in \mathcal{S}.$ It is simple to rewrite $\sup_{t\leq T}|A^N(\phi)_t|$ to see that it is bounded above by a constant dependent on $\phi$ multiplied by $L_T^N/\sqrt{N}$, for $L_t^N$ as defined above Lemma \ref{Quadratictight}. Therefore, by Lemma \ref{Quadratictight}, $\sup_{t\leq T}|A^N(\phi)_t|$ converges uniformly to zero.
\end{proof}

A final error term in $\Upsilon^N$ is given by
$Z^N\in D([0,T],\mathcal{S'})$, defined by $Z^N(t,\phi)=\int_0^t c_{s^-}^N\textrm{d}(\phi',\mu^N_{s})$.
\begin{lemma}
\label{L^Ndif}
    $(Z^N,c^N)$ converges in distribution under the product topology to $(Z,c)\in D([0,T],\mathcal{S}')\times \mathcal{C}([0,T])$. $Z$ is almost surely continuous, and almost surely $Z(t,\phi)=\int_0^t c_{s}\textrm{d}(\phi',\tilde V^\Lambda_{s})$ for all $t\leq T$.
\end{lemma}

\begin{proof}[Proof of Lemma \ref{L^Ndif}]
By analogous argument to \cite[Proposition 5.2]{Distributionalconvergence}, it is sufficient to check that for all finite collections $(\phi_i)_{i=1}^n\in \mathcal{S}$, $((Z^N,\phi_i)_{i=1}^n,c^N)$ converges in distribution under the uniform topology to the process $((\int_0^. c_{s}\textrm{d}(\phi'_i,\tilde V^\Lambda_{s}))_{i=1}^n,c_.)$. We show the result for $n=1$. It is simple to extend to finite $n$.

We apply Ito's formula to $(\phi',\mu^N_{s})$ to obtain that
$$\int_0^t c_{s^-}^Nd(\phi',\mu^N_{s})=-\int_0^tc_{s^-}^N(\phi''',\mu_{s^-}^N)\textrm{d}\Lambda^N_s+\int_0^t\frac{1}{2}c_{s^-}^N(\phi'',\mu_{s^-}^N)\textrm{d}s-\int_0^tc_{s^-}^N\phi'(0)\textrm{d}\Lambda^N_s+E^N_t.$$
$E^N_t$ is an error term which can be shown to converge uniformly to zero, by the same reasoning as for $A^N_t(\phi)$ given in Lemma \ref{anlemma}. Applying the results Theorem \ref{clt} and Lemma \ref{particlesystemconvergencewithoutscaling} to each of the terms in this equation, we deduce that $(\int_0^t c_{s^-}^Nd(\phi',\mu^N_{s}),c^N)$ converges in distribution, under the uniform topology, to $-\int_0^tc_{s}(\phi''',\tilde V^\Lambda_{s^-})\textrm{d}\Lambda_s+\int_0^t\frac{1}{2}c_{s}(\phi'',\tilde V^\Lambda_{s^-})\textrm{d}s-\int_0^tc_s\phi'(0)\textrm{d}\Lambda_s.$ This is identical to $\int_0^tc_sd(\phi',\tilde V^\Lambda_s)$, and is continuous. 
\end{proof}

\begin{proof}[Proof of Lemma \ref{tightnessUpsilon^N}]

\underline{Step 1:}
 We first show that
$$\limsup_{k\rightarrow \infty}\limsup_{N\rightarrow \infty}\mathbb{P}(\sup_{t\leq T}|(\Upsilon^N_t,\phi)|>k)=0.$$

It is clear that this result holds upon replacement of $\Upsilon^N$ by $\tilde \Upsilon^N$. Hence it is sufficient to show that 
$$\limsup_{k\rightarrow \infty}\limsup_{N\rightarrow \infty}\mathbb{P}(\sup_{t\leq T}|(\Upsilon^N_t-\tilde \Upsilon_t^N,\phi)|>k)=0.$$

On the event that $\sup_{s\leq T}|c^N_s|\leq C$,
\begin{align*}
    &\left|(\Upsilon^N_t-\tilde \Upsilon_t^N,\phi)\right|=\left|\frac{1}{\sqrt{N}}\sum_{i=1}^\infty \left(\phi(X_t^{i,N})\mathbbm{1}_{\{\tau_i^N>t\}}-\phi(\tilde X_t^{i,N})\mathbbm{1}_{\{\tilde \tau_i^N>t\}}\right)\right|,\\
    &\leq C||\phi'||_\infty +\left|\frac{1}{\sqrt{N}}\sum_{i=1}^\infty \phi(\tilde X_t^{i,N})\left(\mathbbm{1}_{\{\tau_i^N>t\}}-\mathbbm{1}_{\{\tilde \tau_i^N>t\}}\right)\right|,\\
    &\leq C||\phi'||_\infty + \frac{||\phi||_{\infty}}{\sqrt{N}}\left(\sum_{i=1}^\infty \mathbbm{1}_{\{\tau_i^N>t,\tilde \tau_i^N<t\}}+\sum_{i=1}^\infty\mathbbm{1}_{\{\tau_i^N<t,\tilde \tau_i^N>t\}}\right),\\
    &\leq C||\phi'||_\infty+\sqrt{N}{||\phi||_\infty}(\Gamma_N(\Lambda_.+C/\sqrt{N}\mathbbm{1}_{.\geq 0})_t-\Lambda_t)\\
    &\hspace{1.87cm}+\sqrt{N}||\phi||_\infty(\Lambda_t-\Gamma_N(\Lambda_.-C/\sqrt{N}\mathbbm{1}_{.\geq 0})_t),\\
    &\leq C(||\phi'||_\infty+2||\phi||_\infty)+\sqrt{N}{||\phi||_\infty}|\Gamma_N(\Lambda_.+C/\sqrt{N}\mathbbm{1}_{.\geq 0})_t-\Gamma(\Lambda_.+C/\sqrt{N}\mathbbm{1}_{.\geq 0})_t|\\&\hspace{3.8cm}+\sqrt{N}||\phi||_\infty|\Gamma(\Lambda_.-C/\sqrt{N}\mathbbm{1}_{\{.\geq 0\}})_t-\Gamma_N(\Lambda_.-C/\sqrt{N}\mathbbm{1}_{\{.\geq 0\}})_t|,\\
     &\leq C(||\phi'||_\infty+2||\phi||_\infty)+\sup_{t\leq T}\sqrt{N}{||\phi||_\infty}|\Gamma_N(\Lambda_.+C/\sqrt{N}\mathbbm{1}_{.\geq 0})_t-\Gamma(\Lambda_.+C/\sqrt{N}\mathbbm{1}_{.\geq 0})_t|\\&\hspace{3.8cm}+\sup_{t\leq T}\sqrt{N}||\phi||_\infty|\Gamma(\Lambda_.-C/\sqrt{N}\mathbbm{1}_{\{.\geq 0\}\}})_t-\Gamma_N(\Lambda_.-C/\sqrt{N}\mathbbm{1}_{\{.\geq 0})_t|.
\end{align*}
The processes $\sqrt{N}(\Gamma_N(\Lambda_.\pm C/\sqrt{N})_t-\Gamma(\Lambda_.\pm C/\sqrt{N})_t)$ are martingales of the same form as those used in Lemma \ref{lem:aplus}, and the probability that their maximum on $[0,T]$ is greater than a value $x$ may be controlled through identical arguments to those of Lemma \ref{lem:aplus}. It is then clear that 
$$\limsup_{k\rightarrow \infty}\limsup_{N\rightarrow \infty}\mathbb{P}(\sup_{t\leq T}|(\Upsilon^N_t-\tilde \Upsilon_t^N,\phi)|>k,\sup_{t\leq T}|c^N_t|\leq C)=0.$$
Since $\limsup_{C\rightarrow\infty }\limsup_{N\rightarrow\infty }\mathbb{P}(\sup_{s\leq T}|c^N_s|\leq C)=1$, it then follows that
$$\limsup_{k\rightarrow \infty}\limsup_{N\rightarrow \infty}\mathbb{P}(\sup_{t\leq T}|(\Upsilon^N_t-\tilde \Upsilon_t^N,\phi)|>k)=0.$$
From identical arguments the same can be shown upon replacing $\phi$ by $\phi'$ or $\phi''$.

\underline{Step 2:}
To determine tightness, it is now sufficient to check that for all $\epsilon>0$,
$$\limsup_{\delta\downarrow 0}\limsup_{N\rightarrow\infty}\mathbb{P}(\sup_{t\leq T}\sup_{h\leq \delta\wedge( T-t)}|(\Upsilon^N_{t+h},\phi)-(\Upsilon^N_t,\phi)|>\epsilon)=0.$$

Applying Ito's formula, we determine that
\begin{align}
\label{itoappp}
    (\Upsilon^N_t,\phi)&=(\Upsilon^N_0,\phi)+\frac{1}{\sqrt{N}}\int_0^t\sum_{i=1}^\infty \phi'(X_s^{i,N})\mathbbm{1}_{\{\tau_i^N>s\}}\textrm{d}B_s^{i,N}+\int_0^t(\Upsilon^N_s,\frac{1}{2}\phi'')\textrm{d}s\nonumber\\
    &\qquad+\frac{1}{\sqrt{N}}\sum_{i=1}^\infty \sum_{s\leq t}\phi(X_{s^-}^{i,N}-\Delta\Lambda^N_s)\mathbbm{1}_{\{\tau_i^N>s\}}-\phi(X_{s^-}^{i,N})\mathbbm{1}_{\{\tau_i^N\geq s\}}\nonumber\\
    &\qquad+\sqrt{N}\left(\int_0^t(\phi',\tilde V^\Lambda_s)\textrm{d}\Lambda_s+\phi(0)\Lambda_t\right).
   \nonumber \end{align}
This can be further rewritten as
   \begin{align}
   &(\Upsilon^N_0,\phi)+\frac{1}{\sqrt{N}}\int_0^t\sum_{i=1}^\infty \phi'(X_s^{i,N})\mathbbm{1}_{\{\tau_i^N>s\}}\textrm{d}B_s^{i,N}+\int_0^t(\Upsilon^N_s,\frac{1}{2}\phi'')\textrm{d}s\nonumber\\
    &-\int_0^t(\Upsilon^N,\phi')\textrm{d}\Lambda_s-\phi(0)c^N_t+\frac{1}{\sqrt{N}}\sum_{i=1}^\infty \sum_{s\leq t}\phi(X_{s^-}^{i,N}-\Delta\Lambda^N_s)\mathbbm{1}_{\{\tau_i^N>s\}}-\phi(X_{s^-}^{i,N})\mathbbm{1}_{\{\tau_i^N\geq s\}}\nonumber\\
&\qquad+\sqrt{N}\phi(0)\Lambda^N_t+\int_0^t(\phi',\Upsilon^N_s+\sqrt{N}\tilde V^\Lambda_s)\textrm{d}\Lambda_s.\nonumber
\end{align}
Substituting the definition of $Z^N$ and $A^N$ into this expression, it follows that
\begin{align}
    (\Upsilon^N_t,\phi)&=(\Upsilon^N_0,\phi)+S^N(t,\phi)+\int_0^t(\Upsilon^N_s,\frac{1}{2}\phi'')\textrm{d}s\\
    &\qquad-\int_0^t(\Upsilon^N_s,\phi')\textrm{d}\Lambda_s-\phi(0)c^N_t-(\phi',\mu^N_{t})c^N_t+Z^N(t,\phi)+A^N_t(\phi').\nonumber
\end{align}

To show the result of the current lemma, we will now show that each term in the equation \eqref{itoappp} is suitably bounded. First, by Lemma \ref{anlemma}, $\sup_{t\leq T}|A^N_t(\phi)|$ converges in probability to zero. From the workings of Lemma \ref{jointdist}, and since $\tilde S^N_t(\phi)$ converges uniformly in $t$ in probability, it is simple to deduce that $$\limsup_{\delta\downarrow 0}\limsup_{N\rightarrow\infty}\mathbb{P}(\sup_{t\leq T}\sup_{h\leq \delta\wedge( T-t)}|S^N({t+h},\phi)-S^N(t,\phi)|>\epsilon)=0.$$ An identical result may be deduced for $c_t^N\phi(0)$, $c_t^N(\phi',\mu^N_t)$ and $Z^N(t,\phi)$ by applying the results of Theorem \ref{clt}, Lemma \ref{particlesystemconvergencewithoutscaling} and the workings of Lemma \ref{L^Ndif}.

For $h\leq \delta\wedge (T-t)$, we also bound
\begin{align*}
&\left|\int_t^{t+h}(\Upsilon^N_s,\frac{1}{2}\phi'')\textrm{d}s-\int_t^{t+h}(\Upsilon^N_s,\phi')\textrm{d}\Lambda_s\right|\\
&\qquad\leq \sup_{t\leq T}((\delta +\Lambda_{t+\delta}-\Lambda_t))\left(\sup_{s\leq T}|(\Upsilon_s^N,\phi'')|+\sup_{s\leq T}|(\Upsilon^N_s,\phi')|\right).
\end{align*}
From step one, the second term is stochastically bounded over $N$. The first term converges to 0 as $\delta \rightarrow 0$. Therefore, 
$$\limsup_{\delta\downarrow 0}\limsup_{N\rightarrow\infty}\mathbb{P}\left(\sup_{t\leq T}\sup_{h\leq \delta\wedge( T-t)}\left|\int_t^{t+h}(\Upsilon^N_s,\frac{1}{2}\phi'')\textrm{d}s-\int_t^{t+h}(\Upsilon^N_s,\phi')\textrm{d}\Lambda_s\right|>\epsilon\right)=0.$$
We have thus shown that each term in \eqref{itoappp} possesses the appropriate convergence property, hence so does their sum. 

\underline{Step 3:}
We finally need to show that any subsequential limit is almost surely continuous. We have shown that for all $\epsilon>0$,
$$\limsup_{\delta\downarrow 0}\limsup_{N\rightarrow\infty}\mathbb{P}(\sup_{t\leq T}\sup_{h\leq \delta\wedge( T-t)}|(\Upsilon^N_{t+h},\phi)-(\Upsilon^N_t,\phi)|>\epsilon)=0.$$
This is an open set, so that any subsequential limit $F$ then satisfies
$$\limsup_{\delta\downarrow 0}\mathbb{P}(\sup_{t\leq T}\sup_{h\leq \delta\wedge( T-t)}|F_{t+h}-F_t|>\epsilon)=0,$$
and in particular 
$$\mathbb{P}(\sup_{t\leq T}|F_{t}-F_{t-}|>\epsilon)=0.$$
Since $\epsilon>0$ is arbitrary, almost sure continuity of the process $F$ follows.
\end{proof}

\begin{proof}[Proof of Lemma \ref{itoforspde}]
From Lemmas \ref{anlemma}, \ref{L^Ndif}, and from Corollary \ref{tightUpsilon}, we may take a further subsequence on which $(\Upsilon^N,S^N,c^N,\mu^N,A^N,Z^N)$ converges in distribution to processes $(\Upsilon,S,c,\mu,A,Z)$ such that $A=0,$ $\mu=(\tilde V^\Lambda_t(x)dx)_{t\in [0,T]}$ almost surely. As noted in the proof of \cite[Lemma 2.4]{dawsonmeasuresdes}, the space $D([0,T],\mathcal{S}')$ satisfies the assumptions of the version of the Skorokhod Representation Theorem given as \cite[Theorem 2]{Generalskorokhodrep}. This extends to the product space considered above. Applying \cite[Theorem 2]{Generalskorokhodrep} to $(\Upsilon^N,S^N,c^N,\mu^N,A^N,Z^N)$, we deduce the existence of a further subsequence along which we may obtain
a representation under which the processes converge to $(\Upsilon,S,c,\mu,0,Z)$ almost surely. We shall consider what occurs on the almost sure event where $(\Upsilon^N,S^N,c^N,\mu^N,A^N,Z^N)$ converges, omitting the subsequence notation for brevity.

Fix an arbitrary $t\geq 0$ and $\phi\in \mathcal{S}$. We consider the convergence of each term in the formula \eqref{itoappp}. It is immediate that $(\phi,\Upsilon^N_t)$ converges to $(\phi,\Upsilon_t)$, that $S^N(t,\phi)$ converges to $S(t,\phi),$ that $(\phi,\Upsilon^N_0)$ converges to $(\phi,\Upsilon_0)$, that $c_t^N\phi(0)$ converges to $c_t\phi(0)$, and that $c_t^N(\mu_t^N,\phi')$ converges to $c_t(\phi',\mu_t)$. By Lemmas \ref{anlemma} and \ref{L^Ndif}, it follows that $Z^N_t(\phi)$ converges to $\int_0^tc_sd(\tilde V^\Lambda_s,\phi')$, and that $A^N_t(\phi)$ converges to zero.
From the uniform convergence of $(\phi',\Upsilon^N_t)$ and $(\phi'',\Upsilon^N_t)$, it  follows that $\int_0^t(\Upsilon^N_s,\frac{1}{2}\phi'')\textrm{d}s$ converges to $\int_0^t(\Upsilon_s,\frac{1}{2}\phi'')\textrm{d}s$
  and  $-\int_0^t(\Upsilon^N_s,\phi')\textrm{d}\Lambda_s$ converges to $-\int_0^t(\Upsilon_s,\phi')\textrm{d}\Lambda_s$.

Applying the above results to each term in \ref{itoappp} we may therefore observe,
\begin{align*}
        (\Upsilon_t,\phi)=&(\Upsilon_0,\phi)+\frac{1}{2}\int_0^t(\phi'',\Upsilon_s)\textrm{d}s-\int_0^t(\phi',\Upsilon_s)\textrm{d}\Lambda_s\\
        &\qquad-c_t(\phi',\tilde V^\Lambda_t)+Z_t+S(t,\phi)-c_t\phi(0) \\
        &=(\Upsilon_0,\phi)+\frac{1}{2}\int_0^t(\phi'',\Upsilon_s)\textrm{d}s-\int_0^t(\phi',\Upsilon_s)\textrm{d}\Lambda_s\\
        &\qquad-c_t(\phi',\tilde V^\Lambda_t)+\int_0^tc_s(\phi',\partial_s\tilde V^\Lambda_s)\textrm{d}s+S(t,\phi)-c_t\phi(0). 
\end{align*}
Since $\phi$ and $t$ were arbitrary, we are done.
\end{proof}
\begin{proof}[Proof of Lemma \ref{Quadratictight}]
    We first check that $\mathbb{E}((L_t^N)^2)$ is bounded over $N$.
    As argued in the proof of Lemma \ref{isamg}, we may increase $\tilde V^{\Lambda^N}(s^-,.)$ at a jump time to $||g||_\infty$ and bound the size of the jump at this time via comparison to the busy period of an $M/D/\infty$ queue of arrival rate $||g||_\infty$ and service time one, which is a Borel-Tanner random variable. Furthermore, by using this upper bound, we make the size of jumps independent of the previous behaviour of $\Lambda^N$, since the Poisson processes  $P_{N\Lambda^N_{t^-}}$ are independent. Using this upper bound, we compute
    \begin{align}
       \mathbb{E}((L_T^N)^2)&=\mathbb{E}\left(\left(N\sum_{l\leq T}(\Delta \Lambda^N_l)^2\right)^2\right)\nonumber,\\&\leq \mathbb{E}\left(\frac{1}{N^2} \left(\sum_{l\leq T,\Delta \Lambda^N_l>0}M_4\right)+\frac{1}{N^2}\left(\sum_{ l\leq T} M_2\right)^2\right),\nonumber\\
        &\leq \tilde C\mathbb{E}((\Lambda^N_{T})^2)+\frac{\tilde C}{N}\mathbb{E}(\Lambda^N_T).\label{bdonl}
    \end{align}
    We have used $M_i$ to denote the $i$th moment of the Borel-Tanner distribution of parameter $||g||_\infty$, and $\tilde C$ to denote a generic constant. From Theorem \ref{bdtheorem}, we can further deduce that this is bounded by $\tilde C\Lambda_t^2+\tilde C/N$.
    
We now consider arbitrary $\epsilon\in (0,1/2)$, $t\leq T$, and rewrite
\begin{align}
    L_t^N=&L_t^N\mathbbm{1}_{\{\sup_{s\leq T}|c_s^N|>N^\epsilon\}}+L_\epsilon^N\mathbbm{1}_{\{\sup_{s\leq T}|c_s^N|\leq N^\epsilon\}}-\sum_{\epsilon\leq s\leq t}\Delta \Lambda^N_s\mathbbm{1}_{\{\sup_{s\leq T}|c_s^N|>N^\epsilon\}}\label{termstobound}\\
    &+\Lambda^N_t-\Lambda^N_\epsilon+N\sum_{\epsilon\leq s\leq t}(\Delta \Lambda^N_s)(\Delta \Lambda^N_s-1/N)\mathbbm{1}_{\{\sup_{s\leq T}|c_s^N|\leq N^\epsilon\}}.\nonumber
\end{align}
We shall bound each term on the right hand side of the equation \eqref{termstobound}.

The expectation of the second term is bounded above by $M_2(\Lambda_\epsilon+C_\epsilon/\sqrt{N})$ from the bound on $\mathbb{E}(L_\epsilon^N)$. This converges to zero as $N\rightarrow\infty$ and then $\epsilon\rightarrow 0$. The same holds for the fifth term through similar arguments. 

By Cauchy-Schwarz, the expectation of the first term is bounded above by \\$\sqrt{\mathbb{E}((L_t^N)^2)\mathbb{P}(\sup_{s\leq T}|c_s^N|>N^\epsilon)}$. This is further bounded above by $\tilde C\Lambda_t\exp(-\tilde CN^{2\epsilon})$ by Theorem \ref{bdtheorem}, and the previous bound \eqref{bdonl}. Therefore, the first term in the equation \eqref{termstobound} converges to zero as $N\rightarrow\infty$. The third term can be controlled similarly and converges to zero as $N\rightarrow\infty.$ Further, by Theorem \ref{bdtheorem}, the expectation of the fourth term converges to $\Lambda_t$. Therefore, if the expectation of the final term in \eqref{termstobound} converges to zero as $N\rightarrow\infty$ and $\epsilon\downarrow 0$, we will have established that $\mathbb{E}(L_t^N)$ converges to $\Lambda_t$.

To bound the final term, we note that if $\Delta \Lambda^N_t>1/N$, there must be a point of a Poisson process of intensity $\tilde V^{\Lambda^N}(s^-,.)$ in $[0,1/N]$. From the proof of Lemma \ref{tightN}, $ \sup_{\epsilon\leq s\leq T}V^{\Lambda^N}(s^-,.)\leq \tilde C_\epsilon\left(\tilde V^{\Lambda}(s,.)+\sup_{s\leq T}|c_s^N|/\sqrt{N}\right)$. Therefore, for $\epsilon<s<T$, the conditional probability that there is a point of the Poisson process of intensity $\tilde V^{\Lambda^N}(s^-,.)$ in $[0,1/N]$, given that $\sup_{s\leq T}|c_s^N|\leq N^\epsilon$, is less than 
$\tilde C_\epsilon N^{\epsilon-1/2}+\sup_{\epsilon\leq s\leq T}||\partial_xV(s,.)||_\infty/N$. Further, given that $\Delta \Lambda^N_s>1/N$, the size of $\Delta \Lambda^N_s-2/N$, can then be bounded above by comparison to a Borel-Tanner random variable, in the same manner as for $\Delta \Lambda^N_s$.

Therefore,
\begin{align*}
    &\mathbb{E}\left(N\sum_{\epsilon\leq s\leq t}(\Delta \Lambda^N_s)(\Delta \Lambda^N_s-1/N)\mathbbm{1}_{\{\sup_{s\leq T}|c_s^N|\leq N^\epsilon\}}\right)\\
    &\qquad\leq \frac{1}{N}\mathbb{E}\left(\sum_{\epsilon\leq s\leq t,\Delta \Lambda^N_s>0}\left(M_2+2M_1\right)\left(\tilde C_\epsilon N^{\epsilon-1/2}+\sup_{\epsilon\leq s\leq T}||\partial_xV(s,.)||_\infty/N\right)\right)\\
    &\qquad\leq C_\epsilon \mathbb{E}(\Lambda^N_t)N^{\epsilon-1/2}.
\end{align*}
Sending $N\rightarrow\infty$, this converges to zero. This completes the proof.
\end{proof}
\begin{proof}[Proof of Lemma \ref{lemmahalfwaytodist}]
We first consider $g$ to be a probability density function. We will later extend this to positive measurable functions. Consider a fixed $t>0.$ 
We compute
    \begin{align}
    \label{alabel}
    &N\mathbb{E}\left(\left((\tilde \Lambda^N_t-\Lambda_t)-(M^N_t+\bar M^N_t)\right)^2\right)\nonumber\\
        &\qquad =N\mathbb{E}\left(\left((\tilde \Lambda^N_t-\Lambda_t)-(\Lambda^N_t-\Gamma(\Lambda^N)_t)\right)^2\right),\nonumber\\
        &\qquad=N\mathbb{E}([\tilde \Lambda^N-\Lambda]_t)+N\mathbb{E}([\Lambda^N-\Gamma(\Lambda^N)]_t)-2N\mathbb{E}(\tilde \Lambda^N_t(\Lambda^N_t-\Gamma(\Lambda^N)_t)),\nonumber\\
        &\qquad=\mathbb{E}(\tilde \Lambda^N_t)+N\mathbb{E}\left(\sum_{s\leq t}\left(\Delta \Lambda^N_s-\tilde V^{\Lambda^N}(s^-,[0,\Delta \Lambda^N_s])\right)^2\right)\\
        &\qquad\quad-2N\mathbb{E}\left([H^N,\Lambda^N-\Gamma(\Lambda^N)]_t\right),\nonumber
    \end{align}
where we have defined $H^N_s=\mathbb{E}(\tilde \Lambda^N_t|(\Lambda^N_l)_{l\leq s})$. The second equality follows since both $\tilde \Lambda^N_t-\Lambda_t$ and $\Lambda_t^N-\Gamma(\Lambda^N)_t=M_t^N+\bar M_t^N$ are martingales in their own natural filtrations. The final equality is a simple computation of the quadratic variation and covariation.

The first term in \eqref{alabel}, $\mathbb{E}(\tilde \Lambda^N_t)$, is equal to $\mathbb{E}(\Gamma_N(\Lambda)_t)=\Gamma(\Lambda)_t=\Lambda_t.$ We shall show that the quadratic variation term and the covariation term both converge to $\Lambda_t$, which is sufficient to prove the lemma. 

In order to do this, we first compute $H^N_s$. We label the total number of initial particles as $I^N$. Abusing notation, we shuffle the indices $i=1,...,I^N$ via a uniformly random perturbation so that the particles $(X_t^{i,N})_{i=1,...,I^N}$ are identically distributed, as are $(\tilde X_t^{i,N})_{i=1,...,I^N}$. To consider the behaviour of $\tilde \Lambda^N$ when conditioned on $\Lambda^N$, let us first suppose that $\tau_i^N(\Lambda^N)=l$ for some $i$. Then, $X_0^{i,N}+B_l^{i,N}-\Lambda^N_{l^-}\in [0,\Delta \Lambda^N_l]$ and $X_0^{i,N}+B_x^i-\Lambda^N_x>0$ for $x<l$. Further, conditional on $(\Lambda^N_s)_{s\leq l}$, $(B^{i,N}_{l+x}-B^{i,N}_{l})_{x>0}$ and $X_0^{i,N}+B_l^{i,N}-\Lambda^N_{l^-}$ are independent of $(\Lambda^N_{x}-\Lambda^N_{l})_{x>l}$. This independence follows since if we condition instead on $(\Delta \Lambda^N_s)_{s<l}$ and on $\Delta \Lambda^N_l>0$, the distribution of particles which are not exactly at $0$ forms a Poisson point process, with intensity dependent only on $(\Lambda^N_s)_{s<l}$. Conditioning further on the size of $\Delta \Lambda^N_l$, the set of particles which are above $\Delta \Lambda^N_l$ at the time $l-$ remains Poisson, and is independent of the distribution of the particles below $\Delta \Lambda^N_l$ by the independence property of the Poisson point process.

To simplify the expressions which occur when conditioning we introduce the following notation. First, for a jump time $l$, we denote the distribution of points in $(0,\Delta \Lambda^N_l]$ conditional upon $(\Lambda^N_s)_{s\leq l}$ by $g_l^{N}$. This is identical to conditioning upon $(\Lambda^N_x)_{x>0}$ by the earlier mentioned independence. We label the $k$th dimensional marginal as $g_l^{N,k}$, when $k$ particles are present in $(0,\Delta \Lambda^N_l]$. For a c\'adl\'ag function $f$, a Brownian motion $B$, and $X_0\sim g$ independent of $B$, where $B$ and $X_0$ are both independent of the $\Lambda^N$ particle system, we write $\tau(f)=\inf\{t:X_0+B_t-f_t\leq 0\}$. Finally we write $\mathbb{P}_{N}$ for the probability distribution of $X_0,B$, where we treat $\Lambda^N$ as a deterministic c\'adl\'ag function. 

%hence the later movement of the free boundary is independent of the configuration of points in $[0,\Delta ^N\Lambda_l]$.
We compute
\begin{align}
    &\mathbb{E}(\tilde \Lambda^N_t|(\Lambda^N_l)_{l\leq s})\nonumber\\
    &\qquad=\mathbb{E}\left(\sum_{i\leq I^N}\mathbbm{1}_{\{\tau_i^N(\Lambda)\leq t\}}\mathbbm{1}_{\{\tau_i^N(\Lambda^N)\leq t,X_0^{i,N}+B_{\tau_i^N(\Lambda^N)}-\Lambda^N_{\tau_i^N(\Lambda^N)-}\neq 0\}}\middle|(\Lambda^N_l)_{l\leq s}\right)\nonumber\\
    &\qquad\quad+\mathbb{E}\left(\sum_{i\leq I^N}\mathbbm{1}_{\{\tau_i^N(\Lambda)\leq t\}}\mathbbm{1}_{\{\tau_i^N(\Lambda^N)\leq t,X_0^{i,N}+B_{\tau_i^N(\Lambda^N)}-\Lambda^N_{\tau_i^N(\Lambda^N)-}= 0\}}\middle|(\Lambda^N_l)_{l\leq s}\right)\label{eqcond}\\
    &\qquad\quad+\mathbb{E}\left(\sum_{i\leq I^N}\mathbbm{1}_{\{\tau_i^N(\Lambda)\leq t\}}\mathbbm{1}_{\{\tau_i^N(\Lambda^N)>t\}}\middle | (\Lambda^N)_{l\leq s}\right).\nonumber
    \end{align}
The first term on the right hand side of \eqref{eqcond} is equal to the sum of the probabilities that $\tau(\Lambda)\leq t$, conditional on $\tau(\Lambda^N)=l$, for a jump time $l$ and with the time $l$ distribution of $X_0+B_l-\Lambda_{l-}^N$ being given by $g_l^{N,1}$.
 The first term in \eqref{eqcond} can then be rewritten as 
\begin{align*}
    &\sum_{l\leq s,\Delta \Lambda^N_l>0}(\Delta \Lambda^N_l-1/N)\\
    &\qquad \times \int_0^{\Delta \Lambda^N_l}\mathbb{P}_N(\tau(\Lambda)\leq t|X_0+B_{l}-\Lambda^N_{l^-}=y,\inf_{s<l}(X_0+B_s-\Lambda^N_s)>0)g^{N,1}_{l}(y)\textrm{d}y
    \end{align*}
    The second term on the right hand side of \eqref{eqcond} is equal to the sum of probabilities that $\tau(\Lambda)\leq t$, conditional on $\tau(\Lambda^N)=l$ and on the particle being at $0$ at time $l$. This is given by
\begin{align*}
    \sum_{l\leq s,\Delta \Lambda^N_l>0}\frac{1}{N}\mathbb{P}_N(\tau(\Lambda)\leq t|X_0+B_{l}-\Lambda^N_{l^-}=0,(X_0+B_x-\Lambda^N_x)>0,\space \forall x<l).
    \end{align*}
    The final term on the right hand side of \eqref{eqcond} is $\mathbb{P}_N(\tau(\Lambda)\leq t,\tau(\Lambda^N)>s).$ We therefore write $$H^N_s=\mathbb{E}(\tilde \Lambda^N_t|(\Lambda^N_l)_{l\leq s})=\sum_{l\leq s}P_3^N(t,l)\Delta \Lambda^N_l+\mathbb{P}_N(\tau(\Lambda)\leq t,\tau(\Lambda^N)>s),$$
where $P_3^N(t,l)\in [0,1]$ is defined by
\begin{align*}P_3^N(t,l):=
&\left(1-\frac{1}{N\Delta \Lambda^N_l}\right)\\
&\quad\quad\times\int_0^{\Delta \Lambda^N_l}\mathbb{P}(\tau(\Lambda)\leq t|X_0+B_{l}-\Lambda^N_{l^-}=y,\inf_{x<l}(X_0+B_x-\Lambda^N_x)>0)g^{N,1}_{l}(y)\textrm{d}y
    \\&+\frac{1}{N\Delta \Lambda^N_l}\mathbb{P}(\tau(\Lambda)\leq t|X_0+B_{l}-\Lambda^N_{l^-}=0,(X_0+B_x-\Lambda^N_x)>0\textrm{ for all }x<l).\end{align*}

Using this representation, we can compute the final covariation term in \eqref{alabel} as
\begin{align*}
    &2N\sum_{s\leq t}\left(\left(\Delta \Lambda^N_s-\tilde V^{\Lambda^N}(s^-,[0,\Delta \Lambda^N_s])\right)\left(\Delta \Lambda^N_sP_3^N(t,s)-\mathbb{P}_N(\tau(\Lambda)\leq t,\tau(\Lambda^N)=s\right)
    \right)\\
    &\qquad =2N\sum_{s\leq t}P_3^N(t,s)(\Delta \Lambda^N_s)^2+NO_t^N
    \end{align*}
where $O_t^N$ is a term with $$|O_t^N|\leq 6\sum_{s\leq t}(\Delta \Lambda^N_s)\tilde V^{\Lambda^N}(s^-,[0,\Delta \Lambda^N_s]).$$
From a slight modification to the arguments of Lemma \ref{Quadratictight}, it follows that $\mathbb{E}(N|O^N_t|)\rightarrow 0$ as $N\rightarrow\infty$. Therefore, the limit of the right hand side of \eqref{alabel} is equal to the limit of 
$$\Lambda_t+N\mathbb{E}\left(\sum_{s\leq t}\left(\Delta \Lambda^N_s-\tilde V^{\Lambda^N}(s^-,[0,\Delta \Lambda^N_s])\right)^2\right)-2N\mathbb{E}\left(\sum_{s\leq t}P_3^N(t,s)(\Delta \Lambda^N_s)^2\right).$$
Further, by a slight modification to the proof of Lemma \ref{Quadratictight}, $$\limsup_{N\rightarrow\infty}N\mathbb{E}\left(\sum_{s\leq t}\left(\Delta \Lambda^N_s-\tilde V^{\Lambda^N}(s^-,[0,\Delta \Lambda^N_s])\right)^2\right)=\limsup_{N\rightarrow\infty}\mathbb{E}(L_t^N)=\Lambda_t.$$
Therefore, in order to show that the right hand side of \eqref{alabel} converges to zero it is sufficient to show that
$$\liminf_{N\rightarrow\infty}\mathbb{E}\left(N\sum_{s\leq t}(\Delta \Lambda^N_s)^2P_3^N(t,s))\right)=\Lambda_t,$$
or equivalently that 
$$\limsup_{N\rightarrow\infty}\mathbb{E}\left(N\sum_{s\leq t}(\Delta \Lambda^N_s)^2(1-P_3^N(t,s)))\right)=0.$$

We note that 
$$1-P_3^N(t,s)\leq \mathbb{P}(\inf_{l\leq t-s}c_s^N/\sqrt{N}+B_l>0)\leq (2\sup_{s\leq T}|c_s^N|/\sqrt{N(t-s)} )\wedge 1.$$ 
This follows directly from the definition of $P_3^N(t,s)$, and the strong Markov property of the Brownian motion. Fixing an arbitrary $\epsilon>0$, we then bound
$$\mathbb{E}\left(N\sum_{s\leq t}(\Delta \Lambda^N_s)^2(1-P_3^N(t,s)))\right) \leq \left(\mathbb{E}\left(L^N_\epsilon\right)+\frac{\tilde C_\epsilon}{\sqrt{N}}\sqrt{\mathbb{E}((L_T^N)^2)\mathbb{E}\left(\sup_{s\leq T}|c_s^N|^2\right)}\right).$$

The limsup of the second term on the right hand side is zero from the result of Theorem \ref{bdtheorem} and the arguments in the proof of Lemma \ref{Quadratictight}. The limsup of the first term on the right hand side is $\Lambda_\epsilon$ by Lemma \ref{Quadratictight}. Sending $N\rightarrow\infty$, then $\epsilon\downarrow 0$ we observe that
$$\limsup_{N\rightarrow\infty}\mathbb{E}\left(N\sum_{s\leq t}(\Delta \Lambda^N_s)^2(1-P_3^N(t,s)))\right)=0.$$
In total we have therefore shown that the right hand side of \eqref{alabel} converges to zero, as $N\rightarrow\infty.$

To extend this result to functions $g$ which are not probability densities, we first note that by scaling the result can be extended to any $g\in L^1(\mathbb{R}^+)$. 
To extend to $g\notin L^1(\mathbb{R})$ we then  view only the particles initially below $R$ for $R\in\mathbb{N}$. We set $$\tilde c^{R,N}_t=\sqrt{N}\left(\frac{1}{N}\sum_{i=1}^\infty\mathbbm{1}_{\{\tau_i^N(\Lambda)\leq t,X_0^{i,N}\leq R\}}-\int_0^R\mathbb{P}(\tau^x\leq t)g(x)dx\right)$$ and
$$c^{R,N}_t=\sqrt{N}\left(\frac{1}{N}\sum_{i=1}^\infty\mathbbm{1}_{\{\tau_i^N(\Lambda^N)\leq t,X_0^{i,N}\leq R\}}-\int_0^R\mathbb{P}(\tau^x\leq t)g(x)dx\right).$$
Through an identical proof to Theorem \ref{clt}, we can check that on a subsequence, \\$((c^{R,N},\tilde c^{R,N})_{R\in\mathbb{N}},c^N,\tilde c^N)$ converges in distribution on \\$(D([0,T])^2)^{\oplus\infty}\times D([0,T])\times D([0,T])$ to some $((c^R,\tilde c^R)_{R\in\mathbb{N}},c,\tilde c)$. Further, this limit satisfies that $c_t^R-\tilde c_t^R=\delta \Gamma^R(\Lambda,c^R_t)$ for all $t\leq T$ and all $R\in\mathbb{N}$ almost surely. Here, $\Gamma^R$ is an operator, defined in the same manner as $\Gamma$, but replacing $g$ with $g\mathbbm{1}_{\{.\leq R\}}$. From the definition of the operator $\delta\Gamma(\Lambda,.)$, it follows immediately that if $f_R$ converges to $f$ uniformly,  $\delta\Gamma^R(\Lambda,f_R)\rightarrow\delta\Gamma(\Lambda,f)$ uniformly on compact sets. Further, it is simple to check that $\tilde c^R$ converges to $\tilde c$ as $R\rightarrow\infty$ since the joint distribution is an explicit Gaussian. By considering the form of $c^N-c^{N,R}$, and arguing as in the proof of Theorem \ref{fclthm}, we also determine that $c-c^R$ is a Gaussian process with variance decreasing to zero as $R$ increases to infinity. Therefore $c^R$ converges to $c$ uniformly in probability, hence almost surely on some subsequence. Considering this subsequence, we observe that almost surely $$c_t=\lim_{k\rightarrow\infty }\tilde c_t^{R_k}+\lim_{k\rightarrow\infty} (c^{R_k}_t-\tilde c^{R_k}_t)=\tilde c_t+\lim_{k\rightarrow\infty}\delta \Gamma^{R_k}(\Lambda,c^{R_k})_t=\tilde c_t+\delta\Gamma(\Lambda,c)_t.$$
\end{proof}
\end{document}